\documentclass[11pt]{article}
\usepackage{amssymb}
\usepackage{amsfonts}
\usepackage{amsmath}
\usepackage[OT4]{fontenc}

\newtheorem{thm}{Theorem}[section]
\newtheorem{fact}{Fact}[section]
\newtheorem{lem}{Lemma}[section]
\newtheorem{prop}{Proposition}[section]
\newtheorem{cor}{Corollary}[section]
\newtheorem{rem}{Remark}[section]

\newenvironment{pf*}[1]{\noindent\textit{#1} }{}

\long\def\symbolfootnote[#1]#2{\begingroup%
\def\thefootnote{\fnsymbol{footnote}}\footnote[#1]{#2}\endgroup}

\begin{document}
\begin{center}
{\Large \textbf {Optimal sampling design for global approximation of jump diffusion SDEs 
\symbolfootnote[1]{$\;$ This research was partly supported by the Polish NCN grant - decision No. DEC-2013/09/B/ST1/04275
and by AGH local grant.}
}
}\vspace{0.6cm}

Pawe\l \ Przyby\l owicz

\textit{AGH University of Science
and Technology,\\ 
Faculty of Applied Mathematics,\\
Al. Mickiewicza 30, 30-059 Krakow, Poland,\\
E-mail:} \texttt{pprzybyl@agh.edu.pl}

\end{center}
\vspace{0.2cm}
\begin{abstract}
\noindent
The paper deals with strong global approximation of SDEs driven by two independent processes: a nonhomogeneous Poisson process and a Wiener process. We assume that the jump and diffusion coefficients of the underlying SDE satisfy jump commutativity condition (see Chapter 6.3 in \cite{PBL}). 
We establish the exact convergence rate of minimal errors that can be achieved by arbitrary algorithms based on a finite number of observations of the Poisson and Wiener processes. We consider  classes of methods that use equidistant or nonequidistant sampling of the Poisson and Wiener processes. We provide a construction of optimal methods, based on the classical Milstein scheme, which asymptotically attain the established minimal errors. The analysis implies that methods based on nonequidistant mesh are more efficient than those based on the equidistant mesh.
\end{abstract}
\section{Introduction}
We investigate the global approximation for the following jump diffusion stochastic differential equations (SDEs)
\begin{equation}
	\label{PROBLEM_SDE1}
		\left\{ \begin{array}{ll}
			dX(t)=a(t,X(t))dt+b(t,X(t))dW(t)+c(t,X(t-))dN(t), &t\in [0,T], \\
			X(0)=x_0, 
		\end{array}\right.
\end{equation}
driven by two independent processes: a nonhomogeneous one--dimensional Poisson process $N=\{N(t)\}_{t\in [0,T]}$ with intensity function $\lambda=\lambda(t)>0$ and one-dimensional Wiener process $W$. We assume, without the loss of generality, that $x_0\in\mathbb{R}$. Jump diffusion SDEs (\ref{PROBLEM_SDE1}) appear in various fields such as e.g. physics, biology, engineering and mathematical finance, see, for example, \cite{appl}, \cite{higpla1}, \cite{situ}, \cite{sob} and pages 43-44 in \cite{PBL}. We are interested in efficient algorithms that approximate whole trajectories of $X$ and use only discrete values of the driving Poisson and Wiener processes.

Approximation of stochastic differential equations only driven by a  Wiener process has been widely investigated in the literature.  In that case, upper bounds on the error of defined methods were established, see, for example, \cite{KP}. Lower bounds  were also investigated for the strong approximation in the Wiener (Gaussian) case, see, for example, \cite{hertl}, \cite{HMR3}, \cite{MGHAB}-\cite{AN2} and \cite{PPPM}-\cite{PP6}.

In the jump diffusion case suitable approximation schemes were provided and upper bounds on their errors discussed, for example, in the monograph \cite{PBL} and in the articles  \cite{brupl1}, \cite{gar}, \cite{higpla1}-\cite{higpla3} and \cite{kumsab}. However, according to the author's best knowledge, till now there is only one paper that deals with asymptotic lower bounds and exact rate of convergence of the minimal errors for the global approximation of SDEs with jumps, see \cite{PP7}. In that paper the author considered the pure jumps SDEs (\ref{PROBLEM_SDE1}), i.e., $b\equiv 0$ and $c=c(t)$. We can also mention  \cite{jdpp} where the authors investigated the optimal rate of convergence for the problem of approximating stochastic integrals of regular functions with respect to a homogeneous Poisson process. Here, we extend the approach used in \cite{PP7} in order to cover more general SDEs of the form (\ref{PROBLEM_SDE1}). 

The purpose of this paper is to find lower bounds on the error and to define optimal methods solving (\ref{PROBLEM_SDE1}). In the purely  Gaussian case, similar question were considered, for example,  in  \cite{HMR3} and \cite{MGHAB}. In order to study jump diffusion equations (\ref{PROBLEM_SDE1}) driven by the Poisson and the Wiener processes a new technique is necessary. The main difference, comparing to the Gaussian case, is that we have to use some facts from the theory of stochastic integration with respect to c\`adl\`ag, square integrable martingales, see, for example, \cite{kopp}, \cite{KUO}, \cite{prott} and \cite{situ}. Moreover, we have to face the fact that  when establishing the exact asymptotic constants the intensity of the process $N$ depends on time. This problem does not appear in \cite{PP7}, where the intensity is constant. The another thing is we assume that the coefficients $b$ and $c$ satisfy the  jump commutativity condition. This condition is  widely described and discussed in, for example, Chapter 6.3 in \cite{PBL}. Roughly speaking, it assures  for the construction of the It\^o-Taylor schemes that we do not need to know the exact location of the jump times of the Poisson process $N$. In this paper we widely use this condition when establishing asymptotic lower and upper bounds.

 We consider three classes of approximation schemes denoted by $\chi^{\rm eq}$, $\chi^{\rm noneq*}$ and $\chi^{\rm noneq}$, dependent on the sampling method for trajectories of the processes $N$ and $W$. The class $\chi^{\rm eq}$ contains methods based on the equidistant discretization of $[0,T]$. Methods using the same (but not necessarily equidistant) evaluation points for  $N$ and $W$   belong to a wider class $\chi^{\rm noneq*}$. Methods that can use different, but also not necessarily equidistant, sampling point for the processes $N$ and $W$  belong to  $\chi^{\rm noneq}$. We have  $\chi^{\rm eq}\subset\chi^{\rm noneq*}\subset\chi^{\rm noneq}$.

The main result of the paper, Theorem \ref{OPT_NTH_MIN_ERROR}, states that for fixed $a$, $b$, $c$, $\lambda$, $x_0$ and in the case when the underlying SDE (\ref{PROBLEM_SDE1}) is driven by two processes $N$ and $W$ (i.e., $b\not\equiv 0$ and $c\not\equiv 0$) the following holds
\begin{equation}
\label{INT_1A}
	\limsup\limits_{n\to +\infty} n^{1/2}\cdot \inf\limits_{\bar X_n\in\chi^{\rm noneq}}\Bigl(\mathbb{E}\|X-\bar X_n\|^2_{L^2([0,T])}\Bigr)^{1/2}\leq \frac{1}{\sqrt{6}}\int\limits_0^T\Bigl(\mathbb{E}(\mathcal{Y}(t))\Bigr)^{1/2}dt,
\end{equation}
and 
\begin{equation}
\label{INT_1B}
	\liminf\limits_{n\to +\infty} n^{1/2}\cdot \inf\limits_{\bar X_n\in\chi^{\rm noneq}}\Bigl(\mathbb{E}\|X-\bar X_n\|^2_{L^2([0,T])}\Bigr)^{1/2}\geq \frac{1}{\sqrt{12}}\int\limits_0^T\Bigl(\mathbb{E}(\mathcal{Y}(t))\Bigr)^{1/2}dt,
\end{equation}
where $\mathcal{Y}(t)=|b(t,X(t))|^2+\lambda(t)\cdot\mathbb{E}|c(t,X(t))|^2$, $t\in [0,T]$. In (\ref{INT_1A}) and (\ref{INT_1B}) the method $\bar X_n$ uses at most $n$ evaluations of $N$ and $W$. By taking the infimum we mean that we choose mappings $\{\bar X_n\}_{n\in\mathbb{N}}$ along with discretization points in the best possible way. 
For the subclass $\chi^{\rm noneq*}$ of $\chi^{\rm noneq}$  we have
\begin{equation}
\label{INT_2}
	\lim\limits_{n\to +\infty} n^{1/2}\cdot \inf\limits_{\bar X_n\in\chi^{\rm noneq*}}\Bigl(\mathbb{E}\|X-\bar X_n\|^2_{L^2([0,T])}\Bigr)^{1/2}=\frac{1}{\sqrt{6}}\int\limits_0^T\Bigl(\mathbb{E}(\mathcal{Y}(t))\Bigr)^{1/2}dt,
\end{equation}
while in $\chi^{\rm eq}$  we have that
\begin{equation}
\label{INT_3}
	\lim\limits_{n\to +\infty} n^{1/2}\cdot \inf\limits_{\bar X_n\in\chi^{\rm eq}}\Bigl(\mathbb{E}\|X-\bar X_n\|^2_{L^2([0,T])}\Bigr)^{1/2}=\sqrt{\frac{T}{6}}\Biggl(\int\limits_0^T\mathbb{E}(\mathcal{Y}(t))dt\Biggr)^{1/2}.
\end{equation}
In (\ref{INT_3}) the infimum means that we only choose mappings $\{\bar X_n\}_{n\in\mathbb{N}}$ in the best possible way, while the discretization of $[0,T]$ is fixed and uniform. As we can see, the order of convergence  is $n^{-1/2}$, but the asymptotic constant in (\ref{INT_3}) may be considerably larger than that in (\ref{INT_1A}), (\ref{INT_1B}) and (\ref{INT_2}). In the class $\chi^{\rm noneq}$ we have a small gap between the upper and lower asymptotic constants. We conjecture that the exact rate of convergence of the minimal errors in $\chi^{\rm noneq}$ is the same as for $\chi^{\rm noneq*}$. Note also that if $b\equiv 0$
and $c=c(t)$ then we arrive at results known from \cite{PP7}, while if $c\equiv 0$ and $b\not\equiv 0$ then, for the classes $\chi^{\rm eq}$ and $\chi^{\rm noneq*}$, we restore the results known from \cite{HMR3}, see Remark \ref{rest_know_res}.
\\
The asymptotically optimal scheme is defined by a piecewise linear interpolation of the classical Milstein steps, performed at suitably selected discretization points. 
The discretization points are chosen as quantiles of a distribution corresponding to a density $\psi : [0,T]\to\mathbb{R}_+$. It turns out that in the class $\chi^{\rm noneq*}$ the optimal density $\psi_0$ is proportional to $(\mathbb{E}(\mathcal{Y}(t)))^{1/2}$. The main disadvantage of using such regular sampling is the need of using exact values of quantiles of $(\mathbb{E}(\mathcal{Y}(t)))^{1/2}$ that might be hard to compute in general. In Section 4.1 we present the exact computation of sampling points in the linear case (Merton's model). 

The paper is organized as follows. In Section 2 we give basic notions and definitions. Asymptotic lower bounds on the minimal errors are established in Section 3, while asymptotically optimal methods are defined in Section 4.  We chose such  order of presentation due to the fact that the technique used when proving  the lower bounds in Section 3 suggests definitions of the optimal methods in Section 4. Finally, Appendix contains proofs of auxiliary results used in the paper.
\section{Preliminaries}
Let $T>0$ be a given real number. We denote $\mathbb{N}=\{1,2,\ldots\}$ and $\mathbb{N}_0=\mathbb{N}\cup\{0\}$. Let $(\Omega,\mathcal{F},\mathbb{P})$ be a complete probability space. We consider on it two independent processes: a one-dimensional Wiener process $W=\{W(t)\}_{t\in [0,T]}$ and a one--dimensional nonhomogeneous Poisson process $N=\{N(t)\}_{t\in [0,T]}$ with continuous intensity function $\lambda=\lambda(t)>0$. Let $\{\mathcal{F}_t\}_{t\in [0,T]}$ denote the complete filtration generated by the driving processes $N$ and $W$. We set $\displaystyle{m(t)=\int\limits_0^t\lambda(s)ds}$ and $\Lambda(t,s)=m(t)-m(s)$ for $t,s\in [0,T]$. The process $N$ has independent increments where the  increment $N(t)-N(s)$ has Poisson law with parameter $\Lambda(t,s)$ and $\mathbb{E}(N(t))=m(t)$ for $0\leq s\leq t$, see \cite{GT} or \cite{PBL}. The compensated Poisson process $\tilde N=\{\tilde N(t)\}_{t\in [0,T]}$ is defined as follows
\begin{equation}
	\label{COMP_POISS}
	\tilde N(t)=N(t)-m(t), \quad t\in [0,T],
\end{equation}
which is a zero mean, square integrable $\{\mathcal{F}_t\}_{t\in [0,T]}$-martingale with c\`adl\`ag paths. For a random variable $X:\Omega\to\mathbb{R}$ we write $\|X\|_{L^q(\Omega)}=(\mathbb{E}|X|^q)^{1/q}$, $q\in\{2,4\}$, and $\|X \ | \ \mathcal{G}\|_{L^2(\Omega)}=(\mathbb{E}( |X|^2 \ | \ \mathcal{G} ))^{1/2}$, where $\mathcal{G}$ is a sub-$\sigma$-filed of $\mathcal{F}$.
We say that a continuous function $f:[0,T]\times\mathbb{R}\to \mathbb{R}$ belongs to $C^{0,2}([0,T]\times\mathbb{R})$ if for $j\in\{1,2\}$ the partial derivatives $\partial^j f/\partial y^j=\partial f^j(t,y)/\partial y^j$ exist and are continuous on $(0,T)\times\mathbb{R}$, and can be continuously extended to $[0,T]\times\mathbb{R}$. For a continuous function $f:[0,T]\to\mathbb{R}$ its modulus of continuity is $\bar \omega(f,\delta)=\sup\limits_{t,s\in [0,T], |t-s|\leq\delta}|f(t)-f(s)|$, $\delta\in [0,+\infty)$. If $Y=\{Y(t)\}_{t\in [0,T]}$ is a right-continuous process with left hand limits then we can define $Y(t-):=\lim\limits_{s\to t-}Y(s)$ for all $t\in (0,T]$. We have that $Y(t-)=Y(t)$ if and only if $Y$ is continuous at $t$. For the further properties of c\`adl\`ag mappings used in this paper see, for example, Chapter 2.9 in \cite{appl}.
For $f\in\{b,c\}$  we use the following notation
\begin{eqnarray}
		&& L_1 f(t,y)=b(t,y)\cdot\frac{\partial f}{\partial y}(t,y),\\
		&& L_{-1} f(t,y)=f(t,y+c(t,y))-f(t,y), \quad (t,y)\in [0,T]\times\mathbb{R}.
\end{eqnarray}
We impose the following assumption on the mappings $a:[0,T]\times\mathbb{R}\to\mathbb{R}$,  $b:[0,T]\times\mathbb{R}\to\mathbb{R}$,  $c:[0,T]\times\mathbb{R}\to\mathbb{R}$  and on the intensity function $\lambda:[0,T]\to\mathbb{R}_+$:
\begin{itemize}
	\item [(A)] $f\in C^{0,2}([0,T]\times\mathbb{R})$ for $f\in \{a,b,c\}$.
	\item [(B)] There exists  $K>0$ such that for $f\in\{a,b,c\}$, for all $t,s\in [0,T]$ and all $y,z\in\mathbb{R}$
	\begin{itemize}
		\item [(B1)] $|f(t,y)-f(t,z)|\leq K |y-z|$,
		\item [(B2)] $|f(t,y)-f(s,y)|\leq K(1+|y|)|t-s|$,
		\item [(B3)] $\Bigl|\frac{\partial f}{\partial y}(t,y)-\frac{\partial f}{\partial y}(t,z)\Bigl|\leq K|y-z|$.
	\end{itemize}
	\item [(C)] There exists  $K>0$ such that for $f\in\{b,c\}$, for all $t\in [0,T]$ and all $y,z\in\mathbb{R}$
	$$|L_1f(t,y)-L_1f(t,z)|\leq K|y-z|.$$
	\item [(D)] The diffusion and the jump coefficients satisfy the following \textit{jump commutativity condition} 
	\begin{equation}
		\label{JCC_def}
		L_{-1}b(t,y)=L_1c(t,y),
	\end{equation}
	for all $(t,y)\in [0,T]\times\mathbb{R}$. (We refer to Chapter 6.3 in \cite{PBL} where the condition (\ref{JCC_def}) is widely discussed.)
	\item [(E)] The intensity  $\lambda:[0,T]\to\mathbb{R}_+$ is continuous in $[0,T]$.
\end{itemize}
The  assumptions (B1) and (B2) imply  for $f\in\{a,b,c\}$ and all $(t,y)\in [0,T]\times\mathbb{R}$ that
\begin{equation}
	\label{LIN_GROW1}
	|f(t,y)|\leq K_1 (1+|y|),
\end{equation}
where $K_1>0$ depends only on $f(0,0)$, $K$ and $T$. Moreover, by (B1) and (B3) we have for $f\in\{a,b,c\}$ and all $(t,y)\in [0,T]\times\mathbb{R}$ that
\begin{equation}
	\label{BOUND_AX}
	\Bigl|\frac{\partial^j f}{\partial y^j}(t,y)\Bigl|\leq K, \quad j=1,2.
\end{equation}
From (B1), (\ref{LIN_GROW1}) and (\ref{BOUND_AX}) we get for $f\in\{b,c\}$ and all $(t,y)\in [0,T]\times\mathbb{R}$ that
\begin{equation}
	\label{ling_g_l11}
	\max\{ |L_{-1}f(t,y)|,|L_{1}f(t,y)|\}\leq K_2(1+|y|),
\end{equation}
where $K_2=KK_1$. 
\\
Unless otherwise stated, all unspecified constants appearing in this paper may only
depend on the constant $K$ from the  assumptions (B)-(C), $x_0$, $\|\lambda
\|_{\infty}$, $\|1/\lambda
\|_{\infty}$, $a(0,0)$, $b(0,0)$, $c(0,0)$ and $T$. Moreover, the same symbol might be used to denote different constants.

The assumptions (A)-(E) are rather standard when comparing to those known from the literature concerning approximations of  jump diffusion SDEs, see the comment before Theorem \ref{ERR_MIL_APP}. Only in Section 4.1 we impose additional assumption on the coefficients which, in fact, turns out to be necessary in order to define an optimal sampling from a probabilistic density function.
\\
For $a$, $b$, $c$ and $\lambda$ satisfying (B1), (B2) and (E) the equation (\ref{PROBLEM_SDE1}) has a unique strong solution $X=\{X(t)\}_{t\in [0,T]}$ that is adapted to $\{\mathcal{F}_t\}_{t\in [0,T]}$ and has c\`adl\`ag paths, see \cite{PBL}, \cite{prott} or \cite{sob}. We have also the following  moments estimates for the solution $X$, see, for example, \cite{prott} or \cite{PBL} .
\begin{lem} 
\label{EST_SOL} 
Let us assume that the mappings $a$, $b$, $c$ and $\lambda$ satisfy the assumptions (B1), (B2) and (E). Then there exist  positive constants $C_1$, $C_2$ such that
	\begin{equation}
	\label{EST_SOL_E} 
		\|\sup\limits_{t\in [0,T]}X(t)\|_{L^4(\Omega)}\leq C_1,
	\end{equation}
	and for all $t,s\in [0,T]$
	\begin{equation}
	\label{EST_SOL_CONT} 
		\|X(t)-X(s)\|_{L^2(\Omega)}\leq C_2|t-s|^{1/2}.
	\end{equation}
\end{lem}
The  following result  characterizes the local mean square smoothness of the solution $X$ in the terms of the process  $\mathcal{Y}=\{\mathcal{Y}(t)\}_{t\in [0,T]}$ defined as follows
\begin{equation}
	\label{local_hold_const}
	\mathcal{Y}(t)=|b(t,X(t))|^2+\lambda(t)\cdot|c(t,X(t))|^2, \quad t\in [0,T].
\end{equation}
Of course $\mathcal{Y}$ has c\`adl\`ag paths and it is adapted to $\{\mathcal{F}_{t}\}_{t\in [0,T]}$. (See Fact \ref{exp_cont} in Appendix for the  further properties of $\mathcal{Y}$ used in the paper.)
\begin{prop} 
\label{PROP_REG_SOL} Let us assume that the mappings $a$, $b$, $c$ and $\lambda$ satisfy the assumptions $(B1)$, $(B2)$ and $(E)$. Then for the solution $X$ of (\ref{PROBLEM_SDE1}) we have that for  all $t\in [0,T)$
\begin{equation}
	\label{cond_holder_1}
	\lim\limits_{h\to 0+}\frac{\|X(t+h)-X(t) \ | \ X(t)\|_{L^2(\Omega)}}{h^{1/2}}=(\mathcal{Y}(t))^{1/2},
\end{equation}
almost surely and, in particular, 
\begin{equation}
\label{cond_holder_2}
	\lim\limits_{h\to 0+}\frac{\|X(t+h)-X(t)\|_{L^2(\Omega)}}{h^{1/2}}=\Bigl(\mathbb{E}(\mathcal{Y}(t))\Bigr)^{1/2}.
\end{equation}
\end{prop}
{\bf Proof.} See the Appendix. \ \ \ $\blacksquare$ \\ 
By Proposition \ref{PROP_REG_SOL} the square root of $\mathcal{Y}$ can be interpreted as a conditional H\"older constant of $X$. This local smoothness will reflect in the exact rate of convergence of minimal errors established in Section 4.
A result similar to Proposition \ref{PROP_REG_SOL} for SDEs driven by a multiplicative Wiener process has been obtained in \cite{HMR3}, while for SDEs driven by an additive fractional Brownian motion with the Hurst parameter $H\in (0,1)$ has been shown in Proposition 1 in \cite{AN1}.

The problem considered in the paper is to find an {\it optimal strong global approximation} of the solution $X=\{X(t)\}_{t\in [0,T]}$ of  (\ref{PROBLEM_SDE1}). For any fixed $(a,b,c,\lambda,x_0)$ an approximation of $X=X(a,b,c,\lambda,x_0)$ is given by a method $\bar X=\bar X(a,b,c,\lambda,x_0)$. The method computes the approximation by using some information about the functions $a$, $b$, $c$ and $\lambda$, the Poisson process $N$ and the Wiener process $W$. We consider methods that are based on a finite number of observations of trajectories of the driving processes $N$ and $W$ at suitably chosen points from the interval $[0,T]$. The cost of the method is measured by the total number of evaluations of the processes $N$ and $W$.

We fix $(a,b,c,\lambda,x_0)$ and we consider the corresponding equation (\ref{PROBLEM_SDE1}). Any approximation method $\bar X=\{\bar X_n\}_{n\in\mathbb{N}}$ is defined by three sequences $\bar\varphi=\{\varphi_n\}_{n\in\mathbb{N}}$, $\bar\Delta^Z=\{\Delta^Z_n\}_{n\in\mathbb{N}}$, $Z\in\{N,W\}$, where 
\begin{equation}
	\varphi_n:\mathbb{R}^{2n}\to L^2([0,T]),
\end{equation}
is a measurable mapping and
\begin{equation}
	\label{D_1}
	\Delta^Z_n=\{t^Z_{0,n}, t^Z_{1,n},\ldots ,t^Z_{n,n}\},
\end{equation}
is a partition of $[0,T]$ with
\begin{equation}
	\label{D_2}
		0= t^Z_{0,n}<t^Z_{1,n}<\ldots <t^Z_{n,n}=T,
\end{equation}
for $Z\in\{N,W\}$. We have that $\{0,T\}\subset \Delta^N_n\cap\Delta^W_n$ for all $n$ and, in particular, we might have $\Delta^N_n\cap\Delta^W_n=\emptyset$ for some $n$. The sequences $\bar\Delta^N$, $\bar\Delta^W$ provide (not necessary equidistant)  discretizations of $[0,T]$ used by $N$ and $W$, respectively. Mostly, in the literature, we have that $\bar\Delta^N=\bar\Delta^W$ see, for example, Chapter 6 in \cite{PBL}. Here, mainly for the lower bound, we allow more general discretization.
By 
\begin{equation}
	\label{inf_st_nonad_poiss}
	\mathcal{\bar N}(N,W)=\{\mathcal{N}_{n}(N,W)\}_{n\in\mathbb{N}},
\end{equation}
 we denote a sequence of vectors $\mathcal{N}_n$ of size $2n$, which provides \textit{standard information} with $n$ evaluations of the Poisson process and $n$ evaluations of the Wiener process at the discrete points from $\Delta^N_n\cup\Delta^W_n$, i.e.,
\begin{eqnarray}
	\label{ST_INF_VEC}
	\mathcal{N}_{n}(N,W)&=&\mathcal{N}_{n}(N)\oplus \mathcal{N}_{n}(W)\notag\\
	&:=&[N(t^N_{1,n}), N(t^N_{2,n}),\ldots, N(t^N_{n,n}), W(t^W_{1,n}), W(t^W_{2,n}),\ldots, W(t^W_{n,n})],
\end{eqnarray} 
where $\mathcal{N}_{n}(Z)=[Z(t^Z_{1,n}), Z(t^Z_{2,n}),\ldots, Z(t^Z_{n,n})]$ for $Z\in\{N,W\}$. Recall that $N(0)=W(0)=0$. In particular, the sequences $\bar\varphi$, $\bar\Delta$ may depend on functions $a$, $b$, $c$, $\lambda$ and on $x_0$ but not on trajectories of the processes $N$ and $W$. (Information (\ref{ST_INF_VEC}) uses the same evaluation points for all trajectories of the Poisson and Wiener processes.) Therefore, information (\ref{inf_st_nonad_poiss}) about the processes $N$ and $W$ is {\it nonadaptive}. Moreover, since $\mathcal{N}_{n}(N,W)$ does not have to be contained in $\mathcal{N}_{n+1}(N,W)$, the information (\ref{ST_INF_VEC}) is called {\it nonexpanding}, see \cite{PP5}. We stress that our model of computation covers the regular strong Taylor approximations and it excludes the jump-adapted time discretizations, since we do not assume the knowledge of the  jump times for $N$ (see Chapters 6 and 8 in \cite{PBL}). This restriction reflects our assumption that only nonadaptive standard information is available for the process $N$. 
\\
After computing the information $\mathcal{N}_n(N,W)$, we apply the mapping $\varphi_n$ in order to obtain the $n$th approximation $\bar X_n=\{\bar X_n(t)\}_{t\in [0,T]}$ in the following way
\begin{equation}
	\label{ALG_DEFIN}
	\bar X_n=\varphi_n(\mathcal{N}_{n}(N,W)).
\end{equation}
The $n$th cost of the method $\bar X$ is the total number of evaluations of $N$ and $W$ used by the $n$th approximation $\bar X_n$, defined as follows
\begin{equation}
	\label{NTHCOST}
		cost_n(\bar X)=\left\{ \begin{array}{ll}
			2n, & \hbox{if} \ b\not\equiv 0 \ \hbox{and} \ c\not\equiv 0, \\
			n, & \hbox{if} \ ( b\not\equiv 0 \ \hbox{and} \ c\equiv 0) \ \hbox{or} \ ( b\equiv 0 \ \hbox{and} \ c\not\equiv 0) ,\\
			0, & \hbox{if} \ b\equiv 0 \ \hbox{and} \ c\equiv 0.
		\end{array}\right.
\end{equation}
(If $b\equiv 0$ then we take formally $\mathcal{N}_n(W)$ to be a zero vector and the sequence $\bar\Delta^W$ can be arbitrary; we use analogous convention in the case when $c\equiv 0$.) The set of all methods $\bar X=\{\bar X_n\}_{n\in\mathbb{N}}$, defined as above, is denoted by $\chi^{\rm noneq\it}$. Moreover, we consider the following subclasses of $\chi^{\rm noneq\it}$
\begin{equation}
	\chi^{\rm noneq*\it}=\{\bar X\in\chi^{\rm noneq\it} \ | \ \exists_{n^*_0=n^*_0(\bar X)}:\forall_{n\geq n_0^*} \  \Delta_n^N=\Delta_n^W\},
\end{equation}
and
\begin{equation}
	\chi^{\rm eq\it}=\{\bar X\in\chi^{\rm noneq\it} \ | \ \exists_{n^*_0=n^*_0(\bar X)}:\forall_{n\geq n_0^*} \  \Delta_n^N=\Delta_n^W=\{iT/n \ : \ i=0,1,\ldots,n\}\}.
\end{equation}
Methods based on the sequence of equidistant discretizations (\ref{D_1}) belong to the class $\chi^{\rm eq\it}$ while to the class $\chi^{\rm noneq*\it}$ belong methods that evaluates $N$ and $W$ at the same, possibly nonuniform, sampling points. We have that $\chi^{\rm eq\it}\subset\chi^{\rm noneq*\it}\subset\chi^{\rm noneq\it}$.
\\
The $n$th error of a method $\bar X=\{\bar X_n\}_{n\in\mathbb{N}}$ is defined as 
\begin{equation}
	e_n(\bar X)=\|X-\bar X_n\|_2=\Biggl(\mathbb{E}\int\limits_0^T|X(t)-\bar X_n(t)|^2dt\Biggr)^{1/2}.
\end{equation}
The  $n$th minimal error, in the respective class of methods under consideration, is defined by
\begin{equation}
	\label{N_TH_MINIMAL_ERR}
	e^{\diamond}(n)=\inf_{\bar X\in\chi^{\diamond }}e_n(\bar X), \ \diamond \in\{\rm eq\it, \ \rm noneq*\it, \ \rm noneq\it\}.
\end{equation}
We will investigate the exact rate of convergence of the $n$th minimal errors (\ref{N_TH_MINIMAL_ERR}) together with asymptotic constants. Moreover, we wish to determine  (asymptotically) optimal methods  $\bar X^{\diamond }$, $\diamond \in\{\rm eq\it,\rm noneq*\it ,\rm noneq\it\}$, such that the $n$th errors $e_n(\bar X^{\diamond })$ tend to zero as fast as $e^{\diamond }(n)$ when $n\to +\infty$. 
\section{Asymptotic lower bounds}
In this section we investigate asymptotic lower bounds for the problem (\ref{PROBLEM_SDE1}) in the classes of methods $\chi^\diamond $, $\diamond \in\{\rm eq, \ noneq*, \ noneq\}$. In the next section we give a construction of approximation methods which are asymptotically optimal. Their definitions will be inspired by the technique used for establishing lower bounds given in this section.

We give the definition of the continuous Milstein approximation and we state its properties that we use in order to establish the lower bounds. Moreover, in next section we use it in order to construct asymptotically optimal methods. 
\\
Let $m\in\mathbb{N}$ and
\begin{equation}
\label{MESH_A}
	0 = t_0< t_1 < \ldots < t_m = T ,
\end{equation}	
be an arbitrary discretization of  $[0,T]$. We denote by 
\begin{equation}
	\Delta Z_{i}=Z(t_{i+1})-Z(t_{i}), \quad i=0,1,\ldots,m-1,
\end{equation}	
for $Z\in\{N,W\}$. The continuous Milstein approximation $\tilde X^M_m=\{\tilde X^M_m(t)\}_{t\in [0,T]}$ based on (\ref{MESH_A}) is defined as follows. We denote 
\begin{equation}
	U_i=(t_{i},\tilde X^M_m(t_{i})),
\end{equation}
 and we set
\begin{equation}
	\label{CE_S1}
		\tilde X^M_m(0)=x_0,
\end{equation}
and
\begin{eqnarray}
\label{CE_S2}		
			\tilde X^M_m(t)&=&\tilde X^M_m(t_{i})+a(U_i)\cdot (t-t_{i})+b(U_i)\cdot (W(t)-W(t_{i}))+c(U_i)\cdot (N(t)-N(t_{i}))\notag\\
			&&+L_1b(U_i)\cdot I_{t_i,t}(W,W)+L_{-1}c(U_i)\cdot I_{t_i,t}(N,N)\notag\\
			&&+L_{-1}b(U_i)\cdot I_{t_i,t}(N,W)+L_1c(U_i)\cdot I_{t_i,t}(W,N),
\end{eqnarray}
for $t\in [t_{i},t_{i+1}]$,  $i=0,1,\ldots,m-1$, where
\begin{eqnarray}
	I_{t_i,t}(Y,Z)=\int\limits_{t_i}^{t}\int\limits_{t_i}^{s-} dY(u)dZ(s),
\end{eqnarray}
for $Y,Z\in\{N,W\}$. It is well-known that
\begin{eqnarray}
\label{I_WW}
	&& I_{t_i,t}(W,W)=\frac{1}{2}\Bigl((W(t)-W(t_i))^2-(t-t_i)\Bigr),\\
\label{I_NN}	
	&& I_{t_i,t}(N,N)=\frac{1}{2}\Bigl((N(t)-N(t_i))^2-(N(t)-N(t_i))\Bigr),\\
\label{I_WN2}	
	&& I_{t_i,t}(W,N)=\sum\limits_{k=N(t_i)+1}^{N(t)} W(\tau_k)-W(t_i)\cdot (N(t)-N(t_i)), 
\end{eqnarray}
where $\tau_k$ is the $k$th jump time of $N$, and
\begin{equation}
\label{I_WN}
	I_{t_i,t}(N,W)+I_{t_i,t}(W,N)=(N(t)-N(t_i))\cdot (W(t)-W(t_i)).
\end{equation}
 Moreover, $I_{t_i,t}(W,W)$, $I_{t_i,t}(N,N)$, $I_{t_i,t}(N,W)$, $I_{t_i,t}(W,N)$ are independent of $\mathcal{F}_{t_i}$, see Fact \ref{mul_int_ind} in Appendix.

The main properties of $\tilde X^M_m$ are as follows. For every $m\in\mathbb{N}$ the process $\{\tilde X^M_m(t)\}_{t\in [0,T]}$ is adapted to  $\{\mathcal{F}_t\}_{t\in [0,T]}$ and has c\`adl\`ag paths. The upper bounds on the error of $\tilde X^M_m$ are given in Theorem \ref{ERR_MIL_APP}. Furthermore, under the commutativity condition (\ref{JCC_def}) the random variables $\{\tilde X^M_{m}(t_{i})\}_{i=0}^m$ are measurable with respect to the sigma filed 
\begin{eqnarray*}
	\sigma(\mathcal{N}_{m}(N,W))&=&\sigma(\mathcal{N}_m(N)\oplus \mathcal{N}_m(W))\\
	&:=&\sigma(N(t_1), N(t_2),\ldots, N(t_m), W(t_1), W(t_2),\ldots, W(t_m)).
\end{eqnarray*}
 In particular, this and independence of $N$ and $W$ imply that for all $t\in [t_{i},t_{i+1}]$, $i=0,1,\ldots,m-1$
\begin{eqnarray}
	\label{WEI_PB}
	\tilde X^M_m(t)-\mathbb{E}(\tilde X^M_m(t) \ | \ \mathcal{N}_m(N,W)) &=& b(U_i)\cdot \Bigl( W(t)-\mathbb{E}( W(t) \ | \ \mathcal{N}_m(W))\Bigr)\notag\\
&&+c(U_i)\cdot \Bigl(N(t)-\mathbb{E}( N(t) \ | \ \mathcal{N}_m(N))\Bigr),\notag\\
&&+\tilde R^M_m(t),
\end{eqnarray}
where
\begin{eqnarray}
	\label{def_RM}
	\tilde R^M_m(t)&=&L_1b (U_i)\cdot\Bigl(I_{t_i,t}(W,W)-\mathbb{E}( I_{t_i,t}(W,W) \ | \  \mathcal{N}_m(W))\Bigr)\notag\\
	&&+L_1c (U_i)\cdot\Bigl(I_{t_i,t}(N,W)+I_{t_i,t}(W,N)\notag\\
	&&\quad\quad\quad\quad\quad-\mathbb{E}(I_{t_i,t}(N,W)+I_{t_i,t}(W,N) \ | \  \mathcal{N}_m(N,W))\Bigr)\notag\\
	&&+L_{-1}c (U_i)\cdot\Bigl(I_{t_i,t}(N,N)-\mathbb{E}( I_{t_i,t}(N,N) \ | \  \mathcal{N}_m(N))\Bigr).
\end{eqnarray}
The conditional expectations appearing above can be computed explicitly. Namely, from Lemma 8 in \cite{hertl}  and Lemma \ref{LEM_REGR_N} in Appendix we get by direct calculations
\begin{eqnarray}
\label{COND_I_NN_START}
	&&\mathbb{E}(W(t)-W(t_i) \ | \ \mathcal{N}_m(W))=\Delta W_i\cdot\frac{t-t_{i}}{t_{i+1}-t_{i}},\\
	&&\mathbb{E}(N(t)-N(t_i) \ | \ \mathcal{N}_m(N))=\Delta N_i\cdot\frac{\Lambda(t,t_{i})}{\Lambda(t_{i+1},t_{i})},\\
	&&\mathbb{E}( I_{t_i,t}(W,W) \ | \  \mathcal{N}_m(W))=\Bigl(\frac{t-t_i}{t_{i+1}-t_i}\Bigr)^2\cdot I_{t_i,t_{i+1}}(W,W),\\
	&&\mathbb{E}( I_{t_i,t}(N,W)+I_{t_i,t}(W,N) \ | \  \mathcal{N}_m(N,W))=\Delta N_i\cdot\frac{\Lambda(t,t_i)}{\Lambda(t_{i+1},t_i)}\cdot\Delta W_i\cdot\frac{t-t_i}{t_{i+1}-t_i},\\
	\label{COND_I_NN}
	&&\mathbb{E}( I_{t_i,t}(N,N) \ | \  \mathcal{N}_m(N))=\Bigl(\frac{\Lambda(t,t_i)}{\Lambda(t_{i+1},t_i)}\Bigr)^2\cdot I_{t_i,t_{i+1}}(N,N).
\end{eqnarray}
We stress that for any $m$ the approximation $\{\tilde X^M_m(t)\}_{t\in [0,T]}$ is not an implementable numerical scheme in our model of computation (even under the commutativity condition (\ref{JCC_def})), since computation of a trajectory of $\tilde X^M_m$ requires complete knowledge of a corresponding trajectories of $N$ and $W$. However, if the condition (\ref{JCC_def}) holds,  by (\ref{I_WW}), (\ref{I_NN}) and (\ref{I_WN}), we can compute values of $\tilde X^M_m$ at the discrete points (\ref{MESH_A}) using only function evaluations of $W$ and $N$ at (\ref{MESH_A}). 

In order to characterize asymptotic lower bounds we define
\begin{eqnarray}
	&& C^{\rm noneq}=\frac{1}{\sqrt{6}}\int\limits_0^T\Bigl(\mathbb{E}(\mathcal{Y}(t))\Bigr)^{1/2}dt,\\
	&& C^{\rm eq}=\sqrt{\frac{T}{6}}\cdot\Biggl(\int\limits_0^T \mathbb{E}(\mathcal{Y}(t)) dt\Biggr)^{1/2},
\end{eqnarray}
where the process $\{\mathcal{Y}(t)\}_{t\in [0,T]}$ is defined in (\ref{local_hold_const}).
We have that
\begin{itemize}
	\item[(i)] $0\leq C^{\rm noneq}\leq  C^{\rm eq}$,
	\item [(ii)] $C^{\rm noneq}=C^{\rm eq}$ iff there exists $\gamma\geq 0$ such that for all $t\in [0,T]$
	\begin{equation}
		\mathbb{E}(\mathcal{Y}(t))=\gamma,
	\end{equation}
	\item [(iii)] $C^{\rm eq}=0$ iff $C^{\rm noneq}=0$ iff $b(t,X(t))=0=c(t,X(t))$ for all $t\in [0,T]$ and almost surely. 
\end{itemize}
We have the following result.
\begin{thm}
\label{LOW_B_CONST}
Let us assume that the mappings $a$, $b$, $c$ and $\lambda$ satisfy the assumptions $(A)$-$(E)$. 
\begin{itemize}
	\item [(i)] Let $\bar X$ be an arbitrary method from $\chi^{\rm noneq}$. Then 
	\begin{equation}
	\label{LOW_B_NONEQ}
		\liminf\limits_{n\to+\infty} \ (cost_n(\bar X))^{1/2}\cdot e_n(\bar X)\geq C^{\rm noneq}.
	\end{equation}
	\item [(ii)] Let $\bar X$ be an arbitrary method from $\chi^{\rm noneq*}$. If $b\not\equiv 0$ and $c\not\equiv 0$ then
	\begin{equation}
	\label{LOW_B_NONEQ2}
		\liminf\limits_{n\to+\infty} \ (cost_n(\bar X))^{1/2}\cdot e_n(\bar X)\geq\sqrt{2}\cdot C^{\rm noneq}.
	\end{equation}
	\item [(iii)] Let $\bar X$ be an arbitrary method from $\chi^{\rm eq}$. 
	If $ b\not\equiv 0 \ \hbox{and} \ c\not\equiv 0$ then
\begin{equation}
			\label{LOW_B_EQ2}
		\liminf\limits_{n\to+\infty} \ (cost_n(\bar X))^{1/2}\cdot e_n(\bar X)\geq\sqrt{2}\cdot C^{\rm eq},
	\end{equation}
	else
		\begin{equation}
			\label{LOW_B_EQ1}
		\liminf\limits_{n\to+\infty} \ (cost_n(\bar X))^{1/2}\cdot e_n(\bar X)\geq C^{\rm eq}.
	\end{equation}
\end{itemize}
\end{thm}
{\bf Proof.} We start by  showing (\ref{LOW_B_NONEQ}) in the case when $b\not\equiv 0$ and $c\not\equiv 0$.
Let $\bar X=\{\bar X_n\}_{n\in\mathbb{N}}\in\chi^{\rm noneq}$ be a method based on an arbitrary sequence of discretizations $\bar\Delta^N=\{\Delta^N_n\}_{n\in\mathbb{N}}$ and $\bar\Delta^W=\{\Delta^W_n\}_{n\in\mathbb{N}}$, where each $\Delta^N_n$ and $\Delta^W_n$ is of the form (\ref{D_1}). Every $\bar X_n$ uses information (\ref{ST_INF_VEC}) about the processes $N$ and $W$. Take any sequence $\{m_n\}_{n\in\mathbb{N}}$ of positive integers such that
\begin{equation}
	\label{POS_SEQ_MN}
	\lim\limits_{n\to +\infty}\frac{n^{1/2}}{m_n}=\lim\limits_{n\to +\infty}\frac{m_n}{n}=0.
\end{equation}
By $\hat\Delta=\{\hat\Delta_n\}_{n\in\mathbb{N}}$ we denote a sequence of discretizations given by $\{\hat\Delta_n\}_{n\in\mathbb{N}}=\{\Delta^N_n\cup\Delta^W_n\cup\Delta_n^{\rm eq\it}\}_{n\in\mathbb{N}}$, where every set $\Delta_n^{\rm eq\it}$ of equidistant points is defined by $\Delta_n^{\rm eq\it}=\{jT/m_n \ | \ j=0,1,\ldots,m_n\}$. Hence, for all $n\in\mathbb{N}$,
\begin{equation}
	\hat\Delta_n=\{\hat t_{0,n}, \hat t_{1,n},\ldots,\hat t_{k_n,n}\},
\end{equation}
with
\begin{equation}
	0=\hat t_{0,n}<\hat t_{1,n}<\ldots<\hat t_{k_n,n}=T, \ \hat t_{j,n}\in\Delta^N_n\cup\Delta^W_n\cup\Delta_n^{\rm eq\it}, \ j=0,1,\ldots,k_n,
\end{equation}
and
\begin{equation}
	\label{EST_KN_UL}
	n\leq k_n\leq 2n+m_n-2.
\end{equation}
Therefore, from (\ref{POS_SEQ_MN}) and (\ref{EST_KN_UL}) we have that
\begin{equation}
	\label{EST_KN_UL_1}
		\liminf\limits_{n\to +\infty}\frac{n}{k_n}\geq\frac{1}{2},
\end{equation}
and, since $\Delta^{\rm eq\it}_n\subset\hat\Delta_n$ for all $n\in\mathbb{N}$, 
\begin{equation}
	\label{DIAM_HAT_D}
	\max\limits_{0\leq i\leq k_n-1}(\hat t_{i+1,n}-\hat t_{i,n})\leq \frac{T}{m_n}.
\end{equation}
We denote by $\mathcal{\hat N}(N,W)=\{\mathcal{\hat N}_{n}(N,W)\}_{n\in\mathbb{N}}$, where each vector $\mathcal{\hat N}_{n}(N,W)$ 
consists of the values of $N$ and $W$ at $\hat\Delta_n$, i.e.,
\begin{eqnarray}
	\mathcal{\hat N}_{n}(N,W)&=&\mathcal{\hat N}_{n}(N)\oplus \mathcal{\hat N}_{n}(W)\notag\\
&=&[N({\hat t_{1,n}}), N({\hat t_{2,n}}),\ldots, N({\hat t_{k_n,n}}), W({\hat t_{1,n}}), W({\hat t_{2,n}}),\ldots, W({\hat t_{k_n,n}})].
\end{eqnarray}
Since $\Delta^N_n\cup\Delta^W_n\subset\hat\Delta_n$ for all $n\in\mathbb{N}$, we have that 
\begin{equation}
	\label{N_SUB_N}
\sigma(N(t^N_{1,n}), N(t^N_{2,n}),\ldots, N(t^N_{n,n}), W(t^W_{1,n}), W(t^W_{2,n}),\ldots, W(t^W_{n,n}))\subset\sigma(\mathcal{\hat N}_{n}(N,W)).
\end{equation}
 Let us denote by $\{\tilde X^M_{k_n}\}_{n\in\mathbb{N}}$ the sequence of continuous Milstein approximations (\ref{CE_S1})-(\ref{CE_S2})  based on the sequence of discretizations $\hat\Delta$ and which use the information $\mathcal{\hat N}(N,W)$ about the processes $N$ and $W$. From Theorem \ref{ERR_MIL_APP} and (\ref{DIAM_HAT_D}) we have that
\begin{equation}
	\label{ERROR_HAT_WP}
	\|X-\tilde X^M_{k_n}\|_2\leq C\cdot m_n^{-1},
\end{equation}
where the positive constant $C$ does not depend on $n$. Moreover, let 
\begin{equation}
	\label{brides}
	\hat Z_{n}(t)=Z(t)-\mathbb{E}(Z(t) \ | \ \mathcal{\hat N}_{n}(Z)),
\end{equation}
 for $Z\in\{N,W\}$ and $t\in [0,T]$. Note that for any $t\in [\hat t_{i,n},\hat t_{i+1,n}]$ the random variable $\hat Z_n(t)$ is a convex combination  of $Z(t)-Z(\hat t_{i,n})$ and $-( Z(\hat t_{i+1,n})- Z(t))$.
Hence, $\hat Z_n(t)$ is independent of $\mathcal{F}_{\hat t_{i,n}}$ for all $t\in [\hat t_{i,n},\hat t_{i+1,n}]$ and the processes $\{\hat N_n(t)\}_{t\in [0,T]}$, $\{\hat W_n(t)\}_{t\in [0,T]}$ are independent.  From (\ref{DIAM_HAT_D}), (\ref{N_SUB_N}), (\ref{ERROR_HAT_WP}), (\ref{WEI_PB}) and Lemma \ref{lem_est_rm} we get
\begin{eqnarray}
	\label{EST_LOW_B_NONEQ_2}
	e_n(\bar X)&\geq& \|\bar X_n-\tilde X^{M}_{k_n}\|_{2}-\|X-\tilde X^{M}_{k_n}\|_{2}\notag\\
	\label{LOW_B_TECH_ALG}
&\geq& \|\tilde X^M_{k_n}-\mathbb{E}(\tilde X^M_{k_n} \ | \ \mathcal{\hat N}_n(N,W))\|_2-C\cdot m_n^{-1}\\
&\geq&\Biggl(\sum_{i=0}^{k_n-1}\int\limits_{\hat t_{i,n}}^{\hat t_{i+1,n}}\mathbb{E}|b(\hat t_{i,n},\tilde X^M_{k_n}(\hat t_{i,n}))\cdot \hat W_n(t)+c(\hat t_{i,n},\tilde X^M_{k_n}(\hat t_{i,n}))\cdot \hat N_n(t)|^2dt\Biggr)^{1/2}\notag\\
&&-\Biggl(\sum_{i=0}^{k_n-1}\int\limits_{\hat t_{i,n}}^{\hat t_{i+1,n}}\mathbb{E}|\tilde R^M_{k_n}(t)|^2dt\Biggr)^{1/2}-C\cdot m_n^{-1}\notag\\
\label{EST_LOW_B_NONEQ_22}
&\geq&\Biggl(\sum_{i=0}^{k_n-1}\int\limits_{\hat t_{i,n}}^{\hat t_{i+1,n}}\mathbb{E}|b(\hat t_{i,n},\tilde X^M_{k_n}(\hat t_{i,n}))|^2\cdot \mathbb{E}|\hat W_n(t)|^2+\mathbb{E}|c(\hat t_{i,n},\tilde X^M_{k_n}(\hat t_{i,n}))|^2\cdot \mathbb{E}|\hat N_n(t)|^2dt\Biggr)^{1/2}\notag\\
&&\quad\quad\quad-C\cdot m_n^{-1}.
\end{eqnarray}
Now, we analyze the asymptotic behavior of the first term in (\ref{EST_LOW_B_NONEQ_22}). From Lemma 8 in \cite{hertl} we have that
\begin{equation}
	\label{lem_8_hert}
	\int\limits_{\hat t_{i,n}}^{\hat t_{i+1,n}}\mathbb{E}|\hat W_n(t)|^2 dt=\frac{1}{6} (\hat t_{i+1,n}-\hat t_{i,n})^2.
\end{equation}
For $i=0,1,\ldots,k_n-1$ and $t\in (\hat t_{i,n}, \hat t_{i+1,n})$ we define
\begin{equation}
	H_{i,n}(t)=\frac{\Lambda(t,\hat t_{i,n})\cdot\Lambda(\hat t_{i+1,n},t)}{(t-\hat t_{i,n})(\hat t_{i+1,n}-t)}.
\end{equation}
Of course $H_{i,n}\in C((\hat t_{i,n}, \hat t_{i+1,n}))$ and it can be continuously extended to $[\hat t_{i,n}, \hat t_{i+1,n}]$, since $H(\hat t_{i,n}+)=\lambda(\hat t_{i,n})\cdot\Lambda(\hat t_{i+1,n},\hat t_{i,n})/(\hat t_{i+1,n}-\hat t_{i,n})$ and $H(\hat t_{i+1,n}-)=\lambda(\hat t_{i+1,n})\cdot\Lambda(\hat t_{i+1,n},\hat t_{i,n})/(\hat t_{i+1,n}-\hat t_{i,n})$ are finite.
Therefore, by Lemma \ref{LEM_REGR_N} and from the mean value theorems we get
\begin{eqnarray}	
\label{int_poiss_br}
	&&\int\limits_{\hat t_{i,n}}^{\hat t_{i+1,n}}\mathbb{E}|\hat N_n(t)|^2 dt=\Lambda(\hat t_{i+1,n},\hat t_{i,n})^{-1}\cdot\int\limits_{\hat t_{i,n}}^{\hat t_{i+1,n}}H_{i,n}(t)\cdot (\hat t_{i+1,n}-t)\cdot (t-\hat t_{i,n})dt\notag\\
	&&=\Lambda(\hat t_{i+1,n},\hat t_{i,n})^{-1}\cdot H_{i,n}(\hat d_{i,n})\cdot\int\limits_{\hat t_{i,n}}^{\hat t_{i+1,n}}(\hat t_{i+1,n}-t)\cdot (t-\hat t_{i,n})dt\notag\\
	&&=\frac{1}{6}\frac{\lambda(\hat\alpha_{i,n})\lambda(\hat\beta_{i,n})}{\lambda(\hat\gamma_{i,n})}\cdot (\hat t_{i+1,n}-\hat t_{i,n})^2\leq \frac{1}{6}\cdot\|1/\lambda\|_{\infty}\cdot\|\lambda\|^2_{\infty}\cdot (\hat t_{i+1,n}-\hat t_{i,n})^2,
\end{eqnarray}
for some $\hat d_{i,n},\hat\alpha_{i,n},\hat\beta_{i,n},\hat\gamma_{i,n}\in [\hat t_{i,n},\hat t_{i+1,n}]$, $i=0,1,\ldots,k_n-1$. Next, for $f\in\{b,c\}$ we have from Theorem \ref{ERR_MIL_APP} that
\begin{eqnarray}
\label{dist_XM_X_L}
	&&\Bigl|\mathbb{E}|f(\hat t_{i,n},\tilde X^M_{k_n}(\hat t_{i,n}))|^2-\mathbb{E}|f(\hat t_{i,n}, X(\hat t_{i,n}))|^2\Bigr|\leq\notag\\
	&&\leq C(1+\sup\limits_{t\in [0,T]}\|\tilde X^M_{k_n}(t)\|_{L^2(\Omega)}+\sup\limits_{t\in [0,T]}\| X(t)\|_{L^2(\Omega)})\cdot\|\tilde X^M_{k_n}(\hat t_{i,n})-X(\hat t_{i,n})\|_{L^2(\Omega)}\notag\\
	&&\leq C_1 \max\limits_{0\leq i\leq k_n-1}(\hat t_{i+1,n}-\hat t_{i,n})\leq \frac{C_1 T}{m_n}.
\end{eqnarray}
Therefore, for $(f,Z)\in \{(b,W),(c,N)\}$ we have by (\ref{lem_8_hert}), (\ref{int_poiss_br}) and (\ref{dist_XM_X_L}) that
\begin{equation}
	\lim\limits_{n\to+\infty} n\cdot\sum_{i=0}^{k_n-1}\int\limits_{\hat t_{i,n}}^{\hat t_{i+1,n}}\Bigl(\mathbb{E}|f(\hat t_{i,n},\tilde X^M_{k_n}(\hat t_{i,n}))|^2-\mathbb{E}|f(\hat t_{i,n}, X(\hat t_{i,n}))|^2\Bigr)\cdot\mathbb{E}|\hat Z_n(t)|^2dt=0.
\end{equation}
This together with the H\"older inequality imply
\begin{eqnarray}
\label{EST_LOW_B_HOLD}
	&&\liminf\limits_{n\to +\infty} n \cdot\sum_{i=0}^{k_n-1}\int\limits_{\hat t_{i,n}}^{\hat t_{i+1,n}} \Bigl(\mathbb{E}|b(\hat t_{i,n},\tilde X^M_{k_n}(\hat t_{i,n}))|^2\cdot\mathbb{E}|\hat W_n(t)|^2+\mathbb{E}|c(\hat t_{i,n},\tilde X^M_{k_n}(\hat t_{i,n}))|^2\cdot\mathbb{E}|\hat N_n(t)|^2\Bigr)dt\notag\\
	&&\geq\liminf\limits_{n\to +\infty} n \cdot\sum_{i=0}^{k_n-1}\int\limits_{\hat t_{i,n}}^{\hat t_{i+1,n}} \Bigl(\mathbb{E}|b(\hat t_{i,n},X(\hat t_{i,n}))|^2\cdot\mathbb{E}|\hat W_n(t)|^2+\mathbb{E}|c(\hat t_{i,n},X(\hat t_{i,n}))|^2\cdot\mathbb{E}|\hat N_n(t)|^2\Bigr)dt\notag\\
	&&\geq\liminf\limits_{n\to +\infty}\frac{n}{6 k_n}\cdot\Biggl(\sum_{i=0}^{k_n-1}\Bigl(\mathbb{E}|b(\hat t_{i,n},X(\hat t_{i,n}))|^2+\mathbb{E}|c(\hat t_{i,n},X(\hat t_{i,n}))|^2\cdot\frac{\lambda(\hat\alpha_{i,n})\lambda(\hat\beta_{i,n})}{\lambda(\hat\gamma_{i,n})}\Bigr)^{1/2}\cdot (\hat t_{i+1,n}-\hat t_{i,n})\Biggr)^{2}.\notag
\end{eqnarray}
We have that
\begin{eqnarray}
	\label{def_sn_1}
	\hat S_n &:=&\sum_{i=0}^{k_n-1}\Bigl(\mathbb{E}|b(\hat t_{i,n},X(\hat t_{i,n}))|^2+\mathbb{E}|c(\hat t_{i,n},X(\hat t_{i,n}))|^2\cdot\frac{\lambda(\hat\alpha_{i,n})\lambda(\hat\beta_{i,n})}{\lambda(\hat\gamma_{i,n})}\Bigr)^{1/2}\cdot (\hat t_{i+1,n}-\hat t_{i,n})\notag\\
	&=&\hat S^*_{n}+\hat R_n^*,
\end{eqnarray}
where
\begin{equation}
\label{def_sns_1}
 \hat S^*_{n}:=\sum_{i=0}^{k_n-1}\Bigl(\mathbb{E}(\mathcal{Y}(\hat t_{i,n}))\Bigr)^{1/2}\cdot (\hat t_{i+1,n}-\hat t_{i,n}),
\end{equation}
and, by (\ref{EST_SOL_E}), 
\begin{eqnarray}
	|\hat R^*_{n}|&=&|\hat S_{n}-\hat S^*_{n}|\leq\sum_{i=0}^{k_n-1} \|c(\hat t_{i,n},X(\hat t_{i,n}))\|_{L^2(\Omega)}\cdot\Bigl|\frac{\lambda(\hat\alpha_{i,n})\lambda(\hat\beta_{i,n})}{\lambda(\hat\gamma_{i,n})}-\lambda(\hat t_{i,n})\Bigr|^{1/2}\cdot (\hat t_{i+1,n}-\hat t_{i,n})\notag\\
	&&\leq C_1\cdot\|1/\lambda\|_{\infty}^{1/2}\cdot\|\lambda\|_{\infty}^{1/2}\cdot\sum_{i=0}^{k_n-1}\Bigl(|\lambda(\hat\beta_{i,n})-\lambda(\hat t_{i,n})|^{1/2}+|\lambda(\hat\alpha_{i,n})-\lambda(\hat\gamma_{i,n})|^{1/2}\Bigr)\cdot (\hat t_{i+1,n}-\hat t_{i,n})\notag\\
	&&\leq 2TC_1\cdot\|1/\lambda\|_{\infty}^{1/2}\cdot\|\lambda\|_{\infty}^{1/2}\cdot(\bar\omega(\lambda,T/m_n))^{1/2},
\end{eqnarray}
since $\displaystyle{|\sqrt{x}-\sqrt{y}|\leq\sqrt{|x-y|}}$ for all $x,y\geq 0$. From the uniform continuity of $\lambda$ we get
\begin{equation}
\label{lim_rns_1}
	\lim\limits_{n\to+\infty} \hat R_n^*=0.
\end{equation}
Hence, by (\ref{def_sn_1}), (\ref{def_sns_1}), (\ref{lim_rns_1}) and Fact \ref{exp_cont} (ii) we have
\begin{equation}
	\label{LIM_HAT_S}
	\lim\limits_{n\to +\infty}\hat S_n=\lim\limits_{n\to +\infty}\hat S_n^*=\int\limits_0^T\Bigl(\mathbb{E}(\mathcal{Y}(t))\Bigr)^{1/2}dt.
\end{equation}
Therefore, by (\ref{POS_SEQ_MN}), (\ref{LOW_B_TECH_ALG}), (\ref{EST_LOW_B_HOLD}) and (\ref{LIM_HAT_S}) we obtain
\begin{equation}
	\liminf\limits_{n\to +\infty} \ (cost_n(\bar X))^{1/2}\cdot e_n(\bar X)	\geq\sqrt{2}\cdot\liminf\limits_{n\to +\infty}\Bigl(\frac{n}{6k_n}\Bigr)^{1/2}\cdot \lim\limits_{n\to +\infty} \hat S_n\geq C^{\rm noneq},
\end{equation}
which ends the proof of (\ref{LOW_B_NONEQ}) in the case when $b\not\equiv 0$ and $c\not\equiv 0$. If $( b\not\equiv 0 \ \hbox{and} \ c\equiv 0)$ or $( b\equiv 0 \ \hbox{and} \ c\not\equiv 0)$ then $cost_n(\bar X)=n$,
\begin{equation}
	\label{EST_KN_UL2}
	n\leq k_n\leq n+m_n-1,
\end{equation}
and
\begin{equation}
	\label{EST_KN_UL_2}
		\lim\limits_{n\to +\infty}\frac{n}{k_n}=1,
\end{equation}
which yield
\begin{equation}
	\liminf\limits_{n\to +\infty} \ (cost_n(\bar X))^{1/2}\cdot e_n(\bar X)\geq\lim\limits_{n\to +\infty}\Bigl(\frac{n}{6k_n}\Bigr)^{1/2}\cdot \lim\limits_{n\to +\infty}\hat S_n\geq C^{\rm noneq}.
\end{equation}
For $b\equiv 0$ and $c\equiv 0$  we obtain trivial lower bound. Finally, if $\bar X\in\chi^{\rm noneq*}$, $b\not\equiv 0$ and $c\not\equiv 0$ then $cost_{n}(\bar X)=2n$ and by (\ref{EST_KN_UL_2}) we get
\begin{equation}
	\liminf\limits_{n\to +\infty} \ (cost_n(\bar X))^{1/2}\cdot e_n(\bar X)\geq\sqrt{2}\cdot\lim\limits_{n\to +\infty}\Bigl(\frac{n}{6k_n}\Bigr)^{1/2}\cdot \lim\limits_{n\to +\infty} \hat S_n\geq\sqrt{2}\cdot C^{\rm noneq},
\end{equation}
which completes the proof of (\ref{LOW_B_NONEQ}). The proofs of (\ref{LOW_B_EQ1}) and (\ref{LOW_B_EQ1}) are straightforward modifications of the proofs of (\ref{LOW_B_NONEQ}) and (\ref{LOW_B_NONEQ2}). Hence, we skip it here.\ \ \ $\blacksquare$
\begin{rem}
\label{REM_LOW_B_ARB}
{\rm Theorem \ref{LOW_B_CONST} gives nontrivial lower bounds only in the case when 
\begin{equation}
\label{NONTR_LB}
	\sup\limits_{t\in [0,T]}\mathbb{E}(\mathcal{Y}(t))>0.
\end{equation}
In this case the presented lower bounds still hold even if we allow for methods to use an arbitrary information about $a$, $b$ and $c$, for example, values of partial derivatives or values of arbitrary linear functionals. If  $\mathbb{E}(\mathcal{Y}(t))=0$ for all $t\in [0,T]$ then (\ref{PROBLEM_SDE1}) becomes (almost surely) deterministic ODE. Then  different lower bounds hold, see, for example, \cite{BK2}.  \ \ \ $\square$}
\end{rem}
\section{Asymptotically optimal methods}
We provide definitions of methods that are asymptotically optimal. The construction is inspired by the technique used for establishing the lower bounds in the previous section. We restrict our consideration to approximation methods based on the regular sequences of discretizations generated by a probability density function $\psi$, see \cite{SAYL}. For the density $\psi$ we assume that
\begin{itemize}
	\item [(P1)] $\psi\in C([0,T])$ and $\psi(t)>0$ for all $t\in [0,T]$.
\end{itemize}
We will use the notation $\bar\Delta_{\psi}=\{\Delta_{\psi,n}\}_{n\in\mathbb{N}}$ for a sequence of discretizations generated by a density $\psi$. The knots
\begin{displaymath}
	\Delta_{\psi,n}=\{t_{0,n}, t_{1,n},\ldots, t_{n,n}\},
\end{displaymath}
of the $n$th discretization are given by
\begin{equation}
	\label{DEST_H_DISC}
		\int\limits_{0}^{t_{i,n}}\psi (s)ds=\frac{i}{n}, \ i=0,1,\ldots,n.
\end{equation}
Hence, by choosing such a density $\psi$ one gets a whole sequence of discretizations $\bar\Delta_{\psi}$. For instance, the sequence of equidistant discretizations is obtained by taking $\psi\equiv 1/T$. Since $0<\|\psi\|_{\infty}^{-1}\leq ||1/\psi||_{\infty}<+\infty$, we have for all $n\in\mathbb{N}$ and $i=0,1,\ldots,n-1$
\begin{equation}
	\label{NORMAL_DISC_H}
	\|\psi\|^{-1}_{\infty} n^{-1}\leq t_{i+1,n}-t_{i,n}\leq \|1/\psi\|_{\infty}n^{-1}.
\end{equation}
We now provide a construction of asymptotically optimal approximation methods. The definition of this method is inspired by the dominating term in the estimation (\ref{LOW_B_TECH_ALG}). Denote by $\tilde X^M_{\psi}=\{\tilde X^M_{\psi,n}\}_{n\in\mathbb{N}}$ the sequence of continuous Milstein approximations (\ref{CE_S1})-(\ref{CE_S2}) based on the sequence of discretizations $\bar\Delta_{\psi}$. For a given density $\psi$, we define the method $\bar X_{\psi}^{cM}=\{\bar X_{\psi,n}^{cM}\}_{n\in\mathbb{N}}$ by
\begin{equation}
	\label{COND_MIL}
		\bar X_{\psi,n}^{cM}(t)=\mathbb{E}( \tilde X^{M}_{\psi,n}(t) \ | \ \mathcal{N}_{\psi,n}(N,W)), \quad t\in [0,T],
\end{equation}
where $\mathcal{N}_{\psi,n}(N,W)$ consists of values of the processes $N$ and $W$ at the points $\Delta_{\psi,n}$. (Hence, we formally take  $\bar\Delta^W=\bar\Delta^N=\bar\Delta_{\psi}$.) We call (\ref{COND_MIL}) \textit{the conditional Milstein method}. We have that $\bar X_{\psi,n}^{cM}\in \chi^{\rm noneq*}$. We present an explicit formula for the algorithm (\ref{COND_MIL}) in order to show that it  has a form that is allowed in our model of computation. By (\ref{JCC_def}), (\ref{CE_S2}), (\ref{COND_MIL}) and (\ref{COND_I_NN_START})-(\ref{COND_I_NN}) each term $\bar X_{\psi,n}^{cM}$  can be written as
\begin{eqnarray}
	\label{COND_MILST_EXPL}
	\bar X_{\psi,n}^{cM}(t)=\tilde X^M_{\psi,n}(t_{i,n})&+&a(t_{i,n},\tilde X^M_{\psi,n}(t_{i,n}))\cdot (t-t_{i,n})\notag\\
&+&b(t_{i,n},\tilde X^M_{\psi,n}(t_{i,n}))\cdot\Delta W_{i,n}\cdot\frac{t-t_{i,n}}{t_{i+1,n}-t_{i,n}}\notag\\
&+&c(t_{i,n},\tilde X^M_{\psi,n}(t_{i,n}))\cdot\Delta N_{i,n}\cdot\frac{\Lambda(t,t_{i,n})}{\Lambda(t_{i+1,n},t_{i,n})}\notag\\
&+&L_1b(t_{i,n},\tilde X^M_{\psi,n}(t_{i,n}))\cdot I_{t_{i,n},t_{i+1,n}}(W,W)\cdot\Bigl(\frac{t-t_{i,n}}{t_{i+1,n}-t_{i,n}}\Bigr)^2\notag\\
&+&L_1c(t_{i,n},\tilde X^M_{\psi,n}(t_{i,n}))\cdot\Delta N_{i,n}\cdot\Delta W_{i,n}\cdot\frac{\Lambda(t,t_{i,n})}{\Lambda(t_{i+1,n},t_{i,n})}\cdot\frac{t-t_{i,n}}{t_{i+1,n}-t_{i,n}}\notag\\
&+&L_{-1}c(t_{i,n},\tilde X^M_{\psi,n}(t_{i,n}))\cdot I_{t_{i,n},t_{i+1,n}}(N,N)\cdot\Bigl(\frac{\Lambda(t,t_{i,n})}{\Lambda(t_{i+1,n},t_{i,n})}\Bigr)^2,
\end{eqnarray}
for $t\in [t_{i,n},t_{i+1,n}]$, $i=0,1,\ldots,n-1$ and $\bar X^{cM}_{\psi,n}(0)=x_0$. Note that $\bar X_{\psi,n}^{cM}$ has continuous trajectories and coincides with $\tilde X^M_{\psi,n}$ at the discretization points. In general, the method $\bar X_{\psi,n}^{cM}$ is not equal to the piecewise linear interpolation $\bar X_{\psi,n}^{Lin-M}$ of the classical Mistein steps, defined as
\begin{equation}
	\bar X_{\psi,n}^{Lin-M}(t)=\frac{\tilde X^M_{\psi,n}(t_{i,n})(t_{i+1,n}-t)+\tilde X^M_{\psi,n}(t_{i+1,n})(t-t_{i,n})}{t_{i+1,n}-t_{i,n}},
\end{equation}\
for $t\in [t_{i,n},t_{i+1,n}]$, $i=0,1,\ldots,n-1$, see Remark \ref{CM_EQ_LIN_M}. However, we use the method $\bar X_{\psi,n}^{cM}$ in order to investigate the error of $\bar X_{\psi,n}^{Lin-M}$ and we show in the sequel that they behave asymptotically in the same way. Moreover, for a fixed discretization $\Delta_{\psi,n}$ the method $\bar X_{\psi,n}^{Lin-M}$ does not evaluates $\Lambda$ and its implementation, at least in the case when $\psi\equiv 1/T$, is straightforward.

In the following theorem we give  the exact convergence rate of the errors for the methods $\bar X_{\psi}^{cM}$ and $\bar X_{\psi}^{Lin-M}$ in the terms of the following asymptotic constant
\begin{equation}
	C_{\psi}=\frac{1}{\sqrt{6}}\Biggl(\int\limits_0^T\frac{\mathbb{E}(\mathcal{Y}(t))}{\psi(t)}dt\Biggr)^{1/2}.
\end{equation}
The strategy of the proof goes as follows. First, we analyze the error of the conditional Milstein method $\bar X_{\psi}^{cM}$.
Due to  its definition  given by the conditional expectation (\ref{COND_MIL}) this can be done by using some estimates already established in the proof of Theorem \ref{LOW_B_CONST}. Then we show that $\bar X_{\psi}^{Lin-M}$ is sufficiently close to $\bar X_{\psi}^{cM}$. This will give us the asymptotic error for the piecewise linear interpolation method $\bar X_{\psi}^{Lin-M}$.
\begin{thm} \label{UPP_BOUND_ALGS_H}
Let us assume that the mappings $a$, $b$, $c$, $\lambda$ and $\psi$ satisfy the assumptions $(A)$-$(E)$ and $(P1)$, and let $\bar X_{\psi}\in\{\bar X^{cM}_{\psi}, \bar X^{Lin-M}_{\psi}\}$. Then if $b\not\equiv 0$ and $c\not\equiv 0$
	\begin{equation}
			\label{UPP_B_ALG_CWP_H}
				\lim_{n\to +\infty} (cost_n(\bar X_{\psi}))^{1/2}\cdot e_n(\bar X_{\psi})=\sqrt{2}\cdot C_{\psi},
	\end{equation}
	else
	\begin{equation}
			\label{UPP_B_ALG_CWP_H2}
				\lim_{n\to +\infty} (cost_n(\bar X_{\psi}))^{1/2}\cdot e_n(\bar X_{\psi})=C_{\psi}.
	\end{equation}
\end{thm}
\bf Proof. \rm From Theorem \ref{ERR_MIL_APP}  and (\ref{NORMAL_DISC_H}) we get
\begin{equation}
	\label{APPROX_BB_CWP_DISC_H_T}
	\Bigl|e_n(\bar X^{cM}_{\psi})-\|\tilde X^M_{\psi,n}-\bar X^{cM}_{\psi,n}\|_2\Bigl|\leq e_n(\tilde X^M_{\psi})\leq C\cdot n^{-1},
\end{equation}
where a constant $C>0$ does not depend on $n$. Moreover, the equality (\ref{DEST_H_DISC}) and the integral mean value theorem yield
\begin{equation}
	\label{MEAN_VAL_DEL_H_XI}
	 \forall i=0,1,\ldots,n-1 \ \exists \xi_{i,n}\in [t_{i,n},t_{i+1,n}]: n\cdot (t_{i+1,n}-t_{i,n})=\frac{1}{\psi(\xi_{i,n})}.
\end{equation}
As in  the proof of Theorem \ref{LOW_B_CONST}, we use the notation
\begin{equation}
	\hat Z_{\psi,n}(t)=Z(t)-\mathbb{E}(Z(t) \ | \ \mathcal{N}_{\psi,n}(Z)), \quad t\in [0,T],
\end{equation}
for $Z\in\{N,W\}$. From (\ref{APPROX_BB_CWP_DISC_H_T}), (\ref{MEAN_VAL_DEL_H_XI}), Lemma \ref{lem_est_rm} and by proceeding analogously as in the proof of Theorem \ref{LOW_B_CONST} we arrive at
\begin{eqnarray}
\label{appr_opt_1}
	&&\lim\limits_{n\to +\infty} n\cdot(e_n(\bar X^{cM}_{\psi}))^2=\lim\limits_{n\to +\infty} n\cdot \|\tilde X^M_{\psi,n}-\bar X^{cM}_{\psi,n}\|^2_2\notag\\
&&=\lim\limits_{n\to +\infty} n\cdot\sum\limits_{i=0}^{n-1}\int\limits_{t_{i,n}}^{t_{i+1,n}}\mathbb{E}|b(t_{i,n},\tilde X^M_{\psi,n}(t_{i,n}))\cdot\hat W_{\psi,n}(t)+c(t_{i,n},\tilde X^M_{\psi,n}(t_{i,n}))\cdot\hat N_{\psi,n}(t)|^2dt\notag\\
&&=\lim\limits_{n\to +\infty} n\cdot\sum\limits_{i=0}^{n-1}\int\limits_{t_{i,n}}^{t_{i+1,n}}\Bigl(\mathbb{E}|b(t_{i,n}, X(t_{i,n}))|^2\cdot\mathbb{E}|\hat W_{\psi,n}(t)|^2+\mathbb{E}|c(t_{i,n}, X(t_{i,n}))|^2\cdot\mathbb{E}|\hat N_{\psi,n}(t)|^2\Bigr)dt\notag\\
&&=\frac{1}{6}\lim\limits_{n\to +\infty} \sum\limits_{i=0}^{n-1}\Bigl(\frac{\mathbb{E}|b(t_{i,n}, X(t_{i,n}))|^2}{\psi(\xi_{i,n})}+\frac{\mathbb{E}|c(t_{i,n}, X(t_{i,n}))|^2}{\psi(\xi_{i,n})}\cdot\frac{\lambda(\alpha_{i,n})\lambda(\beta_{i,n})}{\lambda(\gamma_{i,n})}\Bigr)\cdot (t_{i+1,n}-t_{i,n}),
\end{eqnarray}
for some $d_{i,n},\alpha_{i,n},\beta_{i,n},\gamma_{i,n}\in [t_{i,n}, t_{i+1,n}]$, $i=0,1,\ldots,n-1$. Moreover, we have
\begin{eqnarray}
\label{appr_opt_2}
&&\frac{1}{6}\sum\limits_{i=0}^{n-1}\Bigl(\frac{\mathbb{E}|b(t_{i,n}, X(t_{i,n}))|^2}{\psi(\xi_{i,n})}+\frac{\mathbb{E}|c(t_{i,n}, X(t_{i,n}))|^2}{\psi(\xi_{i,n})}\cdot\frac{\lambda(\alpha_{i,n})\lambda(\beta_{i,n})}{\lambda(\gamma_{i,n})}\Bigr)\cdot (t_{i+1,n}-t_{i,n})\notag\\
&&=\frac{1}{6}\sum\limits_{i=0}^{n-1}\frac{\mathbb{E}(\mathcal{Y}(\xi_{i,n}))}{\psi(\xi_{i,n})}\cdot (t_{i+1,n}-t_{i,n})+\sum_{i=1}^3 A_{n,i},
\end{eqnarray}
where
\begin{eqnarray}
	&&A_{n,1}=\frac{1}{6}\sum\limits_{i=0}^{n-1}\frac{\mathbb{E}|b(t_{i,n}, X(t_{i,n}))|^2-\mathbb{E}|b(\xi_{i,n}, X(\xi_{i,n}))|^2}{\psi(\xi_{i,n})}\cdot (t_{i+1,n}-t_{i,n}),\notag\\
	&&A_{n,2}=\frac{1}{6}\sum\limits_{i=0}^{n-1}\frac{\mathbb{E}|c(t_{i,n}, X(t_{i,n}))|^2-\mathbb{E}|c(\xi_{i,n}, X(\xi_{i,n}))|^2}{\psi(\xi_{i,n})}\cdot\frac{\lambda(\alpha_{i,n})\lambda(\beta_{i,n})}{\lambda(\gamma_{i,n})}\cdot (t_{i+1,n}-t_{i,n}),\notag\\
	&&A_{n,3}=\frac{1}{6}\sum\limits_{i=0}^{n-1}\frac{\mathbb{E}|c(\xi_{i,n}, X(\xi_{i,n}))|^2}{\psi(\xi_{i,n})}\cdot\Bigl(\frac{\lambda(\alpha_{i,n})\lambda(\beta_{i,n})}{\lambda(\gamma_{i,n})}-\lambda(\xi_{i,n})\Bigr)\cdot (t_{i+1,n}-t_{i,n}).\notag
\end{eqnarray}
By Fact \ref{exp_cont} (i), (\ref{EST_SOL_E}) and (\ref{NORMAL_DISC_H}) we get
\begin{equation}
	|A_{n,1}|\leq (1/6)\cdot C\cdot\|1/\psi\|^{3/2}_{\infty}\cdot T\cdot n^{-1/2},
\end{equation}
\begin{equation}
	|A_{n,2}|\leq (1/6)\cdot C\cdot\|1/\psi\|^{3/2}_{\infty}\cdot T\cdot\|1/\lambda\|_{\infty}\cdot\|\lambda\|^2_{\infty}\cdot n^{-1/2},
\end{equation}
and
\begin{equation}
	|A_{n,3}|\leq (T/3)\cdot C\cdot\|1/\psi\|_{\infty}\cdot\|1/\lambda\|_{\infty}\cdot\|\lambda\|_{\infty}\cdot \bar\omega(\lambda,\|1/\psi\|_{\infty}\cdot n^{-1}).
\end{equation}
This and the uniform continuity of $\lambda$ imply
\begin{equation}
\label{appr_opt_3}
	\lim\limits_{n\to +\infty} A_{n,j}=0, \quad j=1,2,3.
\end{equation}
By (\ref{appr_opt_1}), (\ref{appr_opt_2}), (\ref{appr_opt_3}) and Fact \ref{exp_cont} (ii) we obtain
\begin{equation}
	\lim\limits_{n\to +\infty} n\cdot (e_n(\bar X^{cM}_{\psi}))^2
=\frac{1}{6}\lim\limits_{n\to +\infty}\sum\limits_{i=0}^{n-1}\frac{\mathbb{E}(\mathcal{Y}(\xi_{i,n}))}{\psi(\xi_{i,n})}\cdot (t_{i+1,n}-t_{i,n})=C^2_{\psi},
\end{equation}
which ends the proof of (\ref{UPP_B_ALG_CWP_H}) for $\bar X_{\psi}=\bar X^{cM}_{\psi}$.

We now analyze the error of $\bar X_{\psi,n}^{Lin-M}$. Note that
\begin{eqnarray}
	\label{TR_LIN_MILST_EXPL}
	\bar R^{M}_{\psi,n}(t)&=&\bar X_{\psi,n}^{cM}(t)-\bar X_{\psi,n}^{Lin-M}(t)=c(t_{i,n},\tilde X^M_{\psi,n}(t_{i,n}))\cdot\Delta N_{i,n}\cdot\Bigl(\frac{\Lambda(t,t_{i,n})}{\Lambda(t_{i+1,n},t_{i,n})}-\frac{t-t_{i,n}}{t_{i+1,n}-t_{i,n}}\Bigr)\notag\\
	&+& L_1b(t_{i,n},\tilde X^M_{\psi,n}(t_{i,n}))\cdot I_{t_{i,n},t_{i+1,n}}(W,W)\cdot\frac{(t-t_{i,n})\cdot (t-t_{i+1,n})}{(t_{i+1,n}-t_{i,n})^2}\notag\\
	&+&L_1c(t_{i,n},\tilde X^M_{\psi,n}(t_{i,n}))\cdot\Delta N_{i,n}\cdot\Delta W_{i,n}\cdot\frac{t_{i,n}-t}{t_{i+1,n}-t_{i,n}}\cdot\frac{\Lambda(t_{i+1,n},t)}{\Lambda(t_{i+1,n},t_{i,n})}\notag\\
	&+&L_{-1}c(t_{i,n},\tilde X^M_{\psi,n}(t_{i,n}))\cdot I_{t_{i,n},t_{i+1,n}}(N,N)\cdot
\Bigl(\Bigl(\frac{\Lambda(t,t_{i,n})}{\Lambda(t_{i+1,n},t_{i,n})}\Bigr)^2-\frac{t-t_{i,n}}{t_{i+1,n}-t_{i,n}}\Bigr),\notag
\end{eqnarray}
for $t\in [t_{i,n},t_{i+1,n}]$, $i=0,1,\ldots,n-1$.  In addition 
\begin{eqnarray}
	\Bigl|\frac{\Lambda(t,t_{i,n})}{\Lambda(t_{i+1,n},t_{i,n})}-\frac{t-t_{i,n}}{t_{i+1,n}-t_{i,n}}\Bigr|&\leq& C\cdot\Bigl|\lambda(t_{i,n})\cdot(t_{i+1,n}-t_{i,n})\cdot\int\limits_{t_{i,n}}^t(\lambda(s)-\lambda(t_{i,n}))ds\notag\\
&+&\lambda(t_{i,n})\cdot (t-t_{i,n})\cdot\int\limits_{t_{i,n}}^{t_{i+1,n}}(\lambda(t_{i,n})-\lambda(s))ds\Bigl|/(t_{i+1,n}-t_{i,n})^2\notag\\
&\leq& C_1\sup\limits_{t,s\in[t_{i,n},t_{i+1,n}]}|\lambda(t)-\lambda(s)|\leq C_1\bar\omega(\lambda,\|1/\psi\|_{\infty}\cdot n^{-1}),
\end{eqnarray}
for $t\in [t_{i,n},t_{i+1,n}]$, $i=0,1,\ldots,n-1$. For $f\in\{b,c\}$ and $j\in\{-1,1\}$ the random variable $L_jf(t_{i,n},\tilde X^M_{\psi,n}(t_{i,n}))$ is  $\mathcal{F}_{t_{i,n}}$-measurable and the estimate (\ref{EST_INT3}) holds for $U_i:=(t_{i,n},\tilde X^M_{\psi,n}(t_{i,n}))$. Hence, it is independent of $I_{t_{i,n},t_{i+1,n}}(N,N)$, $I_{t_{i,n},t_{i+1,n}}(W,W)$ and $\Delta N_{i,n}\cdot\Delta W_{i,n}$. Therefore, by (\ref{NORMAL_DISC_H}) we have  that
\begin{eqnarray}
\label{EST_TR_R}
 \mathbb{E}|\bar R^M_{\psi,n}(t)|^2 &\leq & 4\cdot\Bigl(\Bigl|\frac{\Lambda(t,t_{i,n})}{\Lambda(t_{i+1,n},t_{i,n})}-\frac{t-t_{i,n}}{t_{i+1,n}-t_{i,n}}\Bigr|^2\cdot\mathbb{E}|c(t_{i,n},\tilde X^M_{\psi,n}(t_{i,n}))|^2\cdot\mathbb{E}|\Delta N_{i,n}|^2\notag\\
 &&+\mathbb{E}|L_1b(t_{i,n},\tilde X^M_{\psi,n}(t_{i,n}))|^2\cdot\mathbb{E}|I_{t_{i,n},t_{i+1,n}}(W,W)|^2\notag\\
 &&+\mathbb{E}|L_1c(t_{i,n},\tilde X^M_{\psi,n}(t_{i,n}))|^2\cdot\mathbb{E}|\Delta N_{i,n}|^2\cdot\mathbb{E}|\Delta W_{i,n}|^2\notag\\
 &&+\mathbb{E}|L_{-1}c(t_{i,n},\tilde X^M_{\psi,n}(t_{i,n}))|^2\cdot\mathbb{E}|I_{t_{i,n},t_{i+1,n}}(N,N)|^2\Bigr)\notag\\
 &\leq& C_3\cdot n^{-1}\cdot(\bar\omega(\lambda,\|1/\psi\|_{\infty}\cdot n^{-1}))^2+ C_4 n^{-2},
\end{eqnarray}
for $t\in [t_{i,n},t_{i+1,n}]$, $i=0,1,\ldots,n-1$. 
Since, from (\ref{EST_TR_R})
\begin{equation}
	\Bigl|e_n(\bar X^{Lin-M}_{\psi})-e_n(\bar X^{cM}_{\psi})\Bigr|\leq\|\bar R^{M}_{\psi,n}\|_2\leq C_4\cdot n^{-1/2}\cdot\bar\omega(\lambda,\|1/\psi\|_{\infty}\cdot n^{-1})+C_5\cdot n^{-1},
\end{equation}
 we obtain (\ref{UPP_B_ALG_CWP_H}) for $\bar X_{\psi}=\bar X^{Lin-M}_{\psi}$. This ends the proof. \ \ \ $\blacksquare$\\ \\
Let us now assume that the following additional assumption is satisfied:
\begin{itemize}
	\item [(P2)] $\inf\limits_{t\in [0,T]}\mathbb{E}(\mathcal{Y}(t))>0$.
\end{itemize}
The  methods $\bar X^{cM}_{\psi}$ and $\bar X^{Lin-M}_{\psi}$ obtain the exact rate of convergence $n^{-1/2}$, with the asymptotic constant $C_{\psi}$ which depends on $\psi$. The best density $\psi_0$, which is unique and minimizes $C_{\psi}$ among all positive mappings $\psi\in C([0,T])$ such that $\displaystyle{\int\limits_0^T \psi(t)dt=1}$, is 
	\begin{equation}
		\label{OPT_DEN_H}
			\displaystyle{\psi_0(t)=\frac{\Bigl(\mathbb{E}(\mathcal{Y}(t))\Bigr)^{1/2}}{\int\limits_0^T\Bigl(\mathbb{E}(\mathcal{Y}(s))\Bigr)^{1/2}ds}}, \quad t\in [0,T].
	\end{equation}
(The minimization property of $\psi_0$ follows from the application of the H\"older inequality.) We stress that $\psi_0$ is strictly positive in $[0,T]$ under the additional assumption (P2). Furthermore,
\begin{equation}
	C_{\psi_0}= C^{\rm noneq}, \quad \hbox{and} \quad C_{1/T}= C^{\rm eq}.
\end{equation}
The following fact characterizes the case when the equidistant sampling is the optimal one.
\begin{fact} Let us assume that the mappings $a$, $b$, $c$ and $\lambda$  satisfy the assumptions $(A)$-$(E)$ and $(P2)$. Then the following assertions are equivalent.
\begin{itemize}
	\item [(i)] $\psi_0\equiv 1/T$.
	\item [(ii)] $\displaystyle{\mathbb{E}(\mathcal{Y}(t))=\frac{1}{T}\int\limits_0^T\Bigl(\mathbb{E}(\mathcal{Y}(s))\Bigr)^{1/2}ds}$ for all $t\in [0,T]$.
	\item [(iii)] $C^{\rm noneq}=C^{\rm eq}>0$.
\end{itemize}
\end{fact}
{\bf Proof.} The assertion can easily be shown by proving the implications $(i)\Rightarrow (ii)\Rightarrow (iii)\Rightarrow (i)$ and we left if for the reader. \ \ \ $\blacksquare$\\ \\
From Theorem \ref{UPP_BOUND_ALGS_H} we directly obtain the following result.
\begin{cor}
	\label{UPP_BOUND_ALGS} 
	Let us assume that the mappings $a$, $b$, $c$ and $\lambda$  satisfy the assumptions $(A)$-$(E)$.
	\begin{itemize}
	\item [(i)] Let us moreover assume that the assumption (P2) is satisfied. If $b\not\equiv 0$ and $c\not\equiv 0$ then for $\bar X_{\psi_0}\in\{\bar X^{cM}_{\psi_0},\bar X^{Lin-M}_{\psi_0}\}$ it holds
	\begin{equation}
		\label{UPP_B_ALG_CWP_HOPT1}
			\lim_{n\to +\infty} (cost_n(\bar X_{\psi_0}))^{1/2}\cdot e_n(\bar X_{{\psi}_0})=\sqrt{2}\cdot C^{\rm noneq},
	\end{equation}
	else
	\begin{equation}
		\label{UPP_B_ALG_CWP_HOPT2}
			\lim_{n\to +\infty} (cost_n(\bar X_{\psi_0}))^{1/2}\cdot e_n(\bar X_{{\psi}_0})=C^{\rm noneq}.
	\end{equation}
	\item [(ii)] Let $\bar X_{1/T}\in\{\bar X^{cM}_{1/T},\bar X^{Lin-M}_{1/T}\}$. If $b\not\equiv 0$ and $c\not\equiv 0$ then  it holds
	\begin{equation}
		\label{UPP_B_ALG_CWP_HOPT3}
			\lim_{n\to +\infty} (cost_n(\bar X_{1/T}))^{1/2}\cdot e_n(\bar X_{1/T})=\sqrt{2}\cdot C^{\rm eq},
	\end{equation}
	else
	\begin{equation}
		\label{UPP_B_ALG_CWP_HOPT4}
			\lim_{n\to +\infty} (cost_n(\bar X_{1/T}))^{1/2}\cdot e_n(\bar X_{1/T})=C^{\rm eq}.
	\end{equation}
	\end{itemize}
\end{cor}
Theorem \ref{LOW_B_CONST} and Corollary \ref{UPP_BOUND_ALGS} imply the main result of the paper.
\begin{thm} 
	\label{OPT_NTH_MIN_ERROR}
	Let us assume that the mappings $a$, $b$, $c$ and $\lambda$  satisfy the assumptions $(A)$-$(E)$.
\begin{itemize}
	\item [(i)] Let us additionally assume that the assumption (P2) is satisfied. If ($b\not\equiv 0$ and $c\equiv 0$) or ($b\equiv 0$ and $c\not\equiv 0$) then
		\begin{equation}
		\label{OPT_RATE_NOEQ}
				\lim\limits_{n\to +\infty} n^{1/2}\cdot e^{\rm noneq\it}(n)=C^{\rm noneq},
		\end{equation}
	and the methods $\bar X^{cM}_{\psi_0},\bar X^{Lin-M}_{\psi_0}$, where $\psi_0$ is defined in (\ref{OPT_DEN_H}), are asymptotically optimal in the class $\chi^{\rm noneq}$. If $b\not\equiv 0$ and $c\not\equiv 0$ then
		\begin{equation}
		\label{OPT_RATE_NOEQ2}
				C^{\rm noneq}/\sqrt{2}\leq\liminf\limits_{n\to +\infty} n^{1/2}\cdot e^{\rm noneq\it}(n)\leq \limsup\limits_{n\to +\infty} n^{1/2}\cdot e^{\rm noneq\it}(n)\leq C^{\rm noneq}.
		\end{equation}
	\item [(ii)] If  the assumption (P2) is satisfied, $b\not\equiv 0$ and $c\not\equiv 0$ then
	\begin{displaymath}
			\lim\limits_{n\to +\infty} n^{1/2}\cdot e^{\rm noneq*\it}(n)=C^{\rm noneq},
			\end{displaymath}
		and the methods $\bar X^{cM}_{\psi_0},\bar X^{Lin-M}_{\psi_0}$ are asymptotically optimal in the class $\chi^{\rm noneq*}$.
	\item [(iii)] We have that
			\begin{displaymath}
				\lim\limits_{n\to +\infty} n^{1/2}\cdot e^{\rm eq\it}(n)= C^{\rm eq},
			\end{displaymath}
			and the methods $\bar X^{cM}_{1/T},\bar X^{Lin-M}_{1/T}$ are asymptotically optimal in the class $\chi^{\rm eq}$.
\end{itemize}
\end{thm}
As we can see the optimal rate of convergence of the minimal errors in the classes $\chi^{\rm eq}$ and $\chi^{\rm noneq*}$ is proportional to $n^{-1/2}$, where $n$ is a total number of evaluations of $N$ and $W$. In the class $\chi^{\rm noneq}$ we have a gap between upper and lower asymptotic constants. We conjecture that (\ref{OPT_RATE_NOEQ}) holds also if $b\not\equiv 0$ and $c\not\equiv 0$.
\\
We end this section with the following remarks.
\begin{rem}\rm 
\label{IMPL_CWP_H}
Theorem \ref{OPT_NTH_MIN_ERROR} implies that the error can be reduced asymptotically by the factor
\begin{displaymath}
	C^{\rm eq}/C^{\rm noneq},
\end{displaymath}
if we use the optimal discretization instead of the equidistant one. However, the optimal density $\psi_0$ and the optimal sampling $\{t_{i,n}\}_{i=0}^n$, defined by
\begin{equation}
	\int\limits_{0}^{t_{i,n}}\Bigl(\mathbb{E}(\mathcal{Y}(t))\Bigr)^{1/2}dt=\frac{i}{n}\int\limits_0^T\Bigl(\mathbb{E}(\mathcal{Y}(t))\Bigr)^{1/2}dt, \quad i=0,1,\ldots,n,
\end{equation}
can be computed explicitly only in particular cases see, for example, Section 4.1. Moreover, the additional assumption (P2) is required. We plan to overwhelm these difficulties in the future work. \ \ \ $\square$
\end{rem}
\begin{rem} \label{rest_know_res}
\rm If $c\equiv 0$, $b\not\equiv 0$ and $T=1$ then, for the classes $\chi^{\rm eq}$ and $\chi^{\rm noneq}$, Theorem \ref{OPT_NTH_MIN_ERROR} restores the results of Theorem 2 (iii) and Proposition 2 from \cite{HMR3} in the Gaussian case, while if $c=c(t)$, $b\equiv 0$ and $\lambda=const$ then we get Theorem 4.2 from \cite{PP7} for the pure jump case with an additive Poisson noise. In addition to this paper, in \cite{PP7} the author established a method based on an adaptive stepsize control that does not depend on the knowledge of $\lambda$. The problem of defining such methods for SDEs of the general type (\ref{PROBLEM_SDE1}) will be the topic of our future work. \ \ \ $\square$
\end{rem}
\begin{rem}\rm 
\label{CM_EQ_LIN_M}
We have	$\bar X^{cM}_{\psi}\equiv\bar X^{Lin-M}_{\psi}$, if $\lambda=const$, $b=b(t)$ and $c=c(t)$. \ \ \ $\square$
\end{rem}
\subsection{Linear case - Merton's jump diffusion model}
Let us consider the following SDE 
\begin{equation}
	\label{MERTON_SDE}
		\left\{ \begin{array}{ll}
			dX(t)=rX(t)dt+\sigma X(t) dW(t)+X(t-)dN(t), &t\in [0,T], \\
			X(0)=x_0>0, 
		\end{array}\right.
\end{equation}
that models the stock price in the Merton's model, see \cite{PBL}. We assume  $\lambda$ to be a constant function and $r\in\mathbb{R}$, $\sigma>0$.
The solution of (\ref{MERTON_SDE}) is
\begin{equation}
	X(t)=x_0 \exp\Bigl((r-\frac{1}{2}\sigma^2)t+\sigma W(t)\Bigr)\cdot 2^{N(t)}.
\end{equation}
We denote  $\gamma=r+\sigma^2/2+3\lambda/2$ and we have
\begin{equation}
	\mathbb{E}(\mathcal{Y}(t))=(\sigma^2+\lambda)\cdot x_0^2\cdot e^{2\gamma t}>0, \quad\hbox{for all}\quad t\in [0,T].
\end{equation}
If $\gamma=0$ then the optimal sampling is the equidistant one and $C^{\rm noneq}=C^{\rm eq}=Tx_0\sqrt{\frac{\sigma^2+\lambda}{6}}$. If $\gamma\neq 0$ then we obtain the following optimal sampling for (\ref{MERTON_SDE})
\begin{equation}
	\label{OPT_LIN_C}
	t_{i,n}=\frac{1}{\gamma}\ln \Biggl(\frac{i}{n}\Bigl(e^{\gamma T}-1\Bigr)+1\Biggr), \quad t=0,1,\ldots,n.
\end{equation}
We have that $t_{i,n}\to iT/n$ for $\gamma\to 0$. Since $C^{\rm noneq}/C^{\rm eq}$ behaves as $\sqrt{2/\gamma T}$ when $\gamma T\to +\infty$, we can gain  by using the nonequidistant mesh.
\section{Conclusions}
We investigated the minimal asymptotic errors for strong global approximation of SDEs driven by the Poisson and Wiener processes. We considered the cases of equidistant and nonequidistant sampling of $N$ and $W$. In both cases, we showed that the minimal error tends to zero like $Cn^{-1/2}$, where $C$ is an average in time of  a local H\"older constant of $X$ and $n$ is the number of evaluations of $N$ and $W$. However, the asymptotic constant $C$ in the case of equidistant sampling can be considerably larger than the asymptotic constant when nonuniform mesh is used. We provided a construction of methods that asymptotically achieve the established minimal errors. 

In this paper, we addressed the case when sampling points for the processes $N$ and $W$ are chosen only in the nonadaptive way with respect to $N$ and $W$. Moreover, we assume that the diffusion and jump coefficients satisfied the jump commutativity condition. For the adaptive sampling and non-commutative case preliminary considerations indicate that the direct application of methods developed in this paper is not possible. Further extension of the presented analysis is needed in that case and we postpone this problem to our future work.
\\
\\{\bf Acknowledgments \rm}\\
Part of this work was done at Banff International Research Station for
Mathematical Innovation and Discovery (BIRS), Alberta, Canada, where the author
participated at the workshop  "Approximation of High-Dimensional Numerical
Problems - Algorithms, Analysis and Applications", Fall 2015. I would like to thank
the Staff of the BIRS for great hospitality.
\section{Appendix}
We use the following version of the It\^o formula for semimartingales with jumps, see, for example, \cite{situ} or \cite{prott}.
\begin{lem} 
\label{ITO_FORMULA} Let us assume that the mappings $a$, $b$, $c$ and $\lambda$ satisfy the assumptions $(B1)$, $(B2)$ and $(E)$.
Let a function $U:\mathbb{R}\to\mathbb{R}$ belongs to $C^2(\mathbb{R})$. Then for the solution $X$ of (\ref{PROBLEM_SDE1}) it holds
\begin{eqnarray*}
	U(X(t))=U(X(0))&+&\int\limits_0^t\Bigl(U'(X(s))\cdot a(s,X(s))+\frac{1}{2}\cdot U''(X(s))\cdot b^2(s,X(s))\Bigr)ds\\
	&+&\int\limits_0^tU'(X(s))\cdot b(s,X(s)) dW(s)\\
	&+&\int\limits_0^t\Bigl(U(X(s-)+c(s,X(s-)))-U(X(s-))\Bigr)dN(s).
\end{eqnarray*}
\end{lem}
The proof of the following fact is straightforward.
\begin{fact} 
	\label{exp_cont} 
Let the mappings $a,b,c$ and $\lambda$ satisfy the assumptions $(B1)$, $(B2)$ and $(E)$. 
\begin{itemize}
	\item [(i)] There exists a constant $C_1>0$ such that for all $f\in\{b,c\}$ and  $t,s\in [0,T]$ we have
	\begin{equation}
		\Bigl|\mathbb{E}|f(t,X(t))|^2-\mathbb{E}|f(s,X(s))|^2\Bigr|\leq C_1|t-s|^{1/2}.
	\end{equation}
	\item [(ii)] The mapping 
		\begin{equation}
			[0,T]\ni t\to\mathbb{E}(\mathcal{Y}(t))\in\mathbb{R}_+\cup\{0\},
		\end{equation}
		is continuous.
	\item [(iii)]	There exists a constant $C_2>0$ such that
	\begin{equation}
		\mathbb{E}(\sup\limits_{t\in [0,T]}\mathcal{Y}(t))\leq C_2.
	\end{equation}
\end{itemize}		
\end{fact}	
\begin{fact}
	\label{mul_int_ind}
	\begin{itemize}
		\item [(i)] There exists $C>0$ such that for all $0\leq s\leq t\leq T$ and $Y,Z\in\{N,W\}$ we have
		\begin{equation}
			\label{EST_ITER_INT11}
				\mathbb{E}|I_{s,t}(Y,Z)|^2\leq C(t-s)^2.
\end{equation}
		\item [(ii)] For all $0\leq s\leq t\leq T$ and $Y,Z\in\{N,W\}$ the stochastic integral $I_{s,t}(Y,Z)$ is independent of $\mathcal{F}_s$.
	\end{itemize}	
\end{fact}
{\bf Proof.} The proof of (i) can be straightforwardly delivered from (\ref{COMP_POISS}), (\ref{I_WN}), the isometry for stochastic integrals driven by martingales and by the independence of $W$ and $N$. Hence, we skip it. 

For the proof of (ii) note that directly from (\ref{I_WW}) and (\ref{I_NN}) we get that $I_{s,t}(Y,Y)$, $Y\in\{N,W\}$, is independent of $\mathcal{F}_s$. So the only case of interest is when $(Y,Z)\in \{(N,W), (W,N)\}$.
\\
Fix $s,t\in [0,T]$, $s\leq t$, and let $\Delta_m=\{\alpha_{0,m},\alpha_{1,m},\ldots,\alpha_{m,m}\}$, $m\in\mathbb{N}$, be a sequence of discretizations of $[s,t]$ such that $s=\alpha_{0,m}<\alpha_{1,m}<\ldots<\alpha_{m,m}=t$ and $\lim\limits_{m\to+\infty}\|\Delta_m\|=0$, where $\|\Delta_m\|=\max\limits_{0\leq i\leq m-1}(\alpha_{i+1,m}-\alpha_{i,m})$. Moreover, let
\begin{equation}
	I_{s,t}^m(N,W)=\sum\limits_{j=0}^{m-1}(N(\alpha_{j,m})-N(s))\cdot(W(\alpha_{j+1,m})-W(\alpha_{j,m})).
\end{equation}
We have that
\begin{equation}
	I_{s,t}(N,W)=\lim\limits_{m\to +\infty}I_{s,t}^m(N,W)\quad\hbox{in}\quad L^2(\Omega).
\end{equation}
Therefore, the sequence $\{I_{s,t}^m(N,W)\}_{m\in\mathbb{N}}$ converges also in probability and, by the independence of the increments of $N$ and $W$,  every random variable $I_{s,t}^m(N,W)$ is independent of $\mathcal{F}_s$. Hence, the limit $I_{s,t}(N,W)$ is also independent of $\mathcal{F}_s$. By (\ref{I_WN}) we get that also $I_{s,t}(W,N)$ is independent of $\mathcal{F}_s$. \ \ \ $\blacksquare$\\  \\
{\bf The proof of Proposition \ref{PROP_REG_SOL}.} By the Markov property of the solution $X$ we have that $\|X(t+h)-X(t) \ | \ X(t)\|_{L^2(\Omega)}=\|X(t+h)-X(t) \ | \ \mathcal{F}_t\|_{L^2(\Omega)}$. For all $t\in [0,T)$ and $h>0$ such that $0\leq t< t+h\leq T$ we have 
\begin{eqnarray}
	\label{DIFF_SIE2}
	X(t+h)-X(t)&=&\int\limits_t^{t+h} \Bigl(a(s,X(s))+\lambda(s) c(s,X(s))\Bigr)ds\notag\\
&&+\int\limits_t^{t+h}b(s,X(s))dW(s)+ \int\limits_t^{t+h}c(s,X(s-))d\tilde N(s).
\end{eqnarray}
From (\ref{LIN_GROW1}) and (E) we obtain that
\begin{equation}
	\mathbb{E}\Biggl(\Biggl|\int\limits_t^{t+h} \Bigl(a(s,X(s))+\lambda(s) c(s,X(s))\Bigr)ds\Biggl|^2 \ \Bigl| \ \mathcal{F}_t\Biggr)\leq C_1\cdot  h^2+C_1\cdot h\cdot \mathbb{E}\Bigl(\int\limits_{t}^{t+h} |X(s)|^2 ds \ | \ \mathcal{F}_t\Bigr),
\end{equation}
almost surely. By Theorem 88 in \cite{situ} we obtain for all $t\in [0,T)$ and almost surely 
\begin{eqnarray}
	&&I_1(W):=\mathbb{E}\Biggl(\Biggl|\int\limits_t^{t+h}b(s,X(s))dW(s)\Biggl|^2 \ \Bigl| \ \mathcal{F}_t\Biggr)=\mathbb{E}\Bigl(\int\limits_t^{t+h} |b(s,X(s))|^2ds \ | \ \mathcal{F}_t\Bigr),\\
	&&I_2(N):=\mathbb{E}\Biggl(\Biggl|\int\limits_t^{t+h}c(s,X(s-))d\tilde N(s)\Biggl|^2 \ \Bigl| \ \mathcal{F}_t\Biggr)=\mathbb{E}\Bigl(\int\limits_t^{t+h} \lambda(s)\cdot|c(s,X(s))|^2ds \ | \ \mathcal{F}_t\Bigr),
\end{eqnarray}
and
\begin{eqnarray}	
	&&\mathbb{E}\Biggl(\int\limits_t^{t+h}b(s,X(s))dW(s)\cdot\int\limits_t^{t+h}c(s,X(s-))d\tilde N(s) \ \Bigl| \ \mathcal{F}_t\Biggr)\notag\\
&&\quad\quad\quad\quad\quad=\mathbb{E}\Biggl(\int\limits_t^{t+h}b(s,X(s))\cdot c(s,X(s-))d\langle W,\tilde N\rangle(s) \ \Bigl| \ \mathcal{F}_t\Biggr)=0,
\end{eqnarray}
since $(W(t)\cdot\tilde N(t),\mathcal{F}_t)_{t\in [0,T]}$ is a martingale.
Therefore, by  Minkowski's inequality for conditional expectations (see \cite{doob}), we have that
\begin{eqnarray}
	\label{INEQ_BS22}
	&&\frac{1}{h^{1/2}}\Bigl(I_1(W)+I_2(N)\Bigr)^{1/2}-\Bigl(C_1\cdot h+C_1\cdot \mathbb{E}\Bigl(\int\limits_{t}^{t+h} |X(s)|^2 ds \ | \ \mathcal{F}_t\Bigr)\Bigr)^{1/2}\notag\\
	&&\quad\quad\quad\quad\quad\quad\leq \frac{\|X(t+h)-X(t) \ | \ \mathcal{F}_t\|_{L^2(\Omega)}}{h^{1/2}}\notag\\
	&&\leq \frac{1}{h^{1/2}}\Bigl(I_1(W)+I_2(N)\Bigr)^{1/2}+\Bigl(C_1\cdot h+C_1\cdot \mathbb{E}\Bigl(\int\limits_{t}^{t+h} |X(s)|^2 ds \ | \ \mathcal{F}_t\Bigr)\Bigr)^{1/2},
\end{eqnarray}
almost surely. From (\ref{EST_SOL_E}), Fact \ref{exp_cont} (iii) and the Lebesgue's dominated convergence theorem for conditional expectations (see \cite{doob}) we have for  all $t\in [0,T)$ and almost surely that
\begin{equation}
	\lim\limits_{h\to 0+}\mathbb{E}\Bigl(\int\limits_t^{t+h} |X(s)|^2ds \ | \ \mathcal{F}_t\Bigr)=0,
\end{equation}
and
\begin{equation}
\label{INEQ_BS32}
	\lim\limits_{h\to 0+}\frac{I_1(W)+I_2(N)}{h}=\lim\limits_{h\to 0+}\mathbb{E}\Bigl(\frac{1}{h}\int\limits_t^{t+h} \mathcal{Y}(s)ds \ | \ \mathcal{F}_t\Bigr)=\mathcal{Y}(t),
\end{equation} 
since $X$ and $\mathcal{Y}$ have c\`adl\`ag paths and  $\mathcal{Y}(t)$ is $\mathcal{F}_t$-measurable. This together with (\ref{INEQ_BS22}) yield (\ref{cond_holder_1}). Now, (\ref{cond_holder_2}) follows from (\ref{cond_holder_1}) and Lebesgue's dominated convergence theorem. \ \ \ $\blacksquare$ 
\begin{lem} 
\label{LEM_REGR_N}
Let $m\in\mathbb{N}$ and let
\begin{equation}
	0 = t_0< t_1 < \ldots < t_m = T, 
\end{equation}	
be an arbitrary discretization of the interval $[0, T ]$ and
\begin{equation}
	\mathcal{N}_m(N)=[N(t_1),N(t_2),\ldots,N(t_m)].
\end{equation}
 Then 
for all $i=0,1,\ldots, m-1$ and $t\in [t_i,t_{i+1}]$
\begin{itemize}
	\item [(i)]
	\begin{equation}
		\label{COND_EXP_N}
		\mathbb{E}(N(t) \ | \ \mathcal{N}_m(N))=\frac{N(t_{i+1})\cdot\Lambda(t,t_i)+N(t_i)\cdot\Lambda(t_{i+1},t)}{\Lambda(t_{i+1},t_i)} ,
	\end{equation}
	almost surely,
	\item [(ii)]
	\begin{equation}
		\label{COND2_EXP_N2}
			\mathbb{E}\Bigl(|N(t)-\mathbb{E}(N(t) \ | \ \mathcal{N}_m(N))|^2 \ \Bigl| \ \mathcal{N}_m(N)\Bigr)=(N(t_{i+1})-N(t_{i}))\cdot\frac{\Lambda(t_{i+1},t)\cdot\Lambda(t,t_i)}{(\Lambda(t_{i+1},t_i))^2},
	\end{equation}
	almost surely and, in particular, 
	\begin{equation}
	\label{COND_EXP_N2}
		\mathbb{E}|N(t)-\mathbb{E}(N(t) \ | \ \mathcal{N}_m(N))|^2=\frac{\Lambda(t_{i+1},t)\cdot\Lambda(t,t_i)}{\Lambda(t_{i+1},t_i)}.
	\end{equation}
\end{itemize}	
\end{lem}
{\bf Proof.} For $t=t_i$, $i=0,1,\ldots, m$, we directly get (\ref{COND_EXP_N}), (\ref{COND2_EXP_N2}) and (\ref{COND_EXP_N2}).
By the results of \cite{bonnual}, from the fact that the process $N$ has independent increments and by direct calculations we obtain that conditioned on $\mathcal{N}_m(N)$ and for $t\in (t_i,t_{i+1})$, $i=0,1,\ldots,m-1$ the increment $N(t)-N(t_i)$ is a binomial random variable with the number of trials $N(t_{i+1})-N(t_{i})$ and with the probability of success in each trial equal to $\displaystyle{\frac{\Lambda(t,t_i)}{\Lambda(t_{i+1},t_{i})}}$. Now, the rest of proof goes analogously as the proof of Lemma 3.1 in \cite{PP7}. \ \ \ $\blacksquare$\\ \\
We provide a result concerning an upper bound on the error for the continuous Milstein approximation $\tilde X^{M}_m$. A similar result has been shown in Theorem 6.4.1 in \cite{PBL}, however, under slightly stronger assumptions. In particular, we do not assume the existence of continuous partial derivative $\partial f/\partial t$ for $f\in \{a,b,c\}$ and we do not impose here any Lipschitz conditions on the second partial derivative of $f=f(t,y)$, $f\in \{a,b,c\}$, with respect to $y$. Moreover, we consider here nonstationary Poisson process, while in \cite{PBL} Theorem 6.4.1 has been proven for stationary  point processes.
\begin{thm} 
\label{ERR_MIL_APP} 
Let us assume that the mappings $a$, $b$, $c$ and $\lambda$ satisfy the assumptions $(A)$, $(B)$, $(C)$ and $(E)$. Let $m\in\mathbb{N}$ and let (\ref{MESH_A}) be an arbitrary discretization of the interval $[0,T]$. Then for the continuous Milstein approximation $\tilde X^M_m$, based on the mesh (\ref{MESH_A}), we have that
\begin{equation}
	\label{EST_EULER_APP1}
	\sup\limits_{t\in [0,T]}\|\tilde X^M_m(t)\|_{L^2(\Omega)}\leq C_1,
\end{equation}
and
\begin{equation}
\label{EST_EULER_APP2}
	\sup\limits_{t\in [0,T]}\|X(t)-\tilde X^M_m(t)\|_{L^2(\Omega)}\leq C_2\cdot\max\limits_{0\leq i\leq m-1}(t_{i+1}-t_i),
\end{equation}
where $C_1,C_2>0$ do not depend on $m$.
\end{thm}
{\bf Proof.}  Recall that  $U_i=(t_i,\tilde X^M_m(t_i))$ and $L_jf(U_i)$ is $\mathcal{F}_{t_i}$-measurable for $f\in\{b,c\}$, $j\in\{-1,1\}$. First, we show that
\begin{equation}
\label{finit_2mom_Mil}
	\sup\limits_{t\in [0,T]}\|\tilde X^M_m(t)\|_{L^2(\Omega)}<+\infty.
\end{equation}
We proceed by induction. Let us assume that $\|\tilde X^M_m(t_k)\|_{L^2(\Omega)}<+\infty$ for $k=0,1,\ldots,i$ and some $i$. (The assumption is fulfilled for $i=0$.)
By (\ref{LIN_GROW1}), (\ref{ling_g_l11}), (\ref{CE_S2}) and Fact \ref{mul_int_ind} we have for all $k=0,1,\ldots,i$ and $t\in [t_k,t_{k+1}]$  that
\begin{equation}
	\|\tilde X^M_m(t)-\tilde X^M_m(t_k)\|_{L^2(\Omega)}\leq C(1+\|\tilde X^M_m(t_k)\|_{L^2(\Omega)})(t-t_k)^{1/2}.
\end{equation}
Hence, $\sup\limits_{t\in [t_i,t_{i+1}]}\|\tilde X^M_m(t)\|_{L^2(\Omega)}<+\infty$ and, in particular, $\|\tilde X^M_m(t_{i+1})\|_{L^2(\Omega)}<+\infty$. Therefore, we get $\max\limits_{i=0,1,\ldots,n}\|\tilde X^M_m(t_{i})\|_{L^2(\Omega)}<+\infty$ and (\ref{finit_2mom_Mil}).
\\
We now justify (\ref{EST_EULER_APP2}). The solution $X=X(t)$ of (\ref{PROBLEM_SDE1}) and the continuous Milstein approximation $\tilde X^M_m=\tilde X^M_m(t)$ can be written as
\begin{eqnarray}
	&& X(t)=x_0+A(t)+B(t)+C(t),\\
	&& \tilde X^M_m(t)=x_0+\tilde A^M_m(t)+\tilde B^M_m(t)+\tilde C^M_m(t),
\end{eqnarray}
where
\begin{eqnarray}
	&&A(t)=\int\limits_0^t\sum\limits_{i=0}^{m-1}a(s,X(s))\cdot\mathbf{1}_{(t_i,t_{i+1}]}(s)ds,\\
	&&B(t)=\int\limits_0^t\sum\limits_{i=0}^{m-1}b(s,X(s))\cdot\mathbf{1}_{(t_i,t_{i+1}]}(s)dW(s),\\
	&&C(t)=\int\limits_0^t\sum\limits_{i=0}^{m-1}c(s,X(s-))\cdot\mathbf{1}_{(t_i,t_{i+1}]}(s)dN(s),
\end{eqnarray}
and
\begin{eqnarray}
	&&\tilde A^M_m(t)=\int\limits_0^t\sum\limits_{i=0}^{m-1}a(U_i)\cdot\mathbf{1}_{(t_i,t_{i+1}]}(s)ds,\notag\\
	&&\tilde B^M_m(t)=\int\limits_0^t\sum\limits_{i=0}^{m-1}\Bigl(b(U_i)+\int\limits_{t_i}^s L_1b(U_i)dW(u)+\int\limits_{t_i}^sL_{-1}b(U_i)dN(u)\Bigr)\cdot\mathbf{1}_{(t_i,t_{i+1}]}(s)dW(s),\notag\\
	&&\tilde C^M_m(t)=\int\limits_0^t\sum\limits_{i=0}^{m-1}\Bigl(c(U_i)+\int\limits_{t_i}^s L_1 c(U_i)dW(u)+\int\limits_{t_i}^{s-}L_{-1} c(U_i)dN(u)\Bigr)\cdot\mathbf{1}_{(t_i,t_{i+1}]}(s)dN(s).\notag
\end{eqnarray}
We have  for all $t\in [0,T]$ that
\begin{equation}
	\label{DIFF_A_AE}
	\mathbb{E}|A(t)-\tilde A^M_m(t)|^2\leq 3\Bigl(\mathbb{E}|\tilde A^M_{m,1}(t)|^2+\mathbb{E}|\tilde A^M_{m,2}(t)|^2+\mathbb{E}|\tilde A^M_{m,3}(t)|^2\Bigr),
\end{equation}
where
\begin{eqnarray}
	&&\mathbb{E}|\tilde A^M_{m,1}(t)|^2=\mathbb{E}\Bigl|\int\limits_0^t\sum\limits_{i=0}^{m-1}\Bigl(a(s,X(s))-a(t_i,X(s))\Bigr)\cdot\mathbf{1}_{(t_i,t_{i+1}]}(s)ds\Bigl|^2,\\
	&&\mathbb{E}|\tilde A^M_{m,2}(t)|^2=\mathbb{E}\Bigl|\int\limits_0^t\sum\limits_{i=0}^{m-1}\Bigl(a(t_i,X(s))-a(t_i,X(t_i))\Bigr)\cdot\mathbf{1}_{(t_i,t_{i+1}]}(s)ds\Bigl|^2,\\
	\label{EST_AE_3}
	&&\mathbb{E}|\tilde A^M_{m,3}(t)|^2=\mathbb{E}\Bigl|\int\limits_0^t\sum\limits_{i=0}^{m-1}\Bigl(a(t_i,X(t_i))-a(t_i,\tilde X^M_m(t_i))\Bigr)\cdot\mathbf{1}_{(t_i,t_{i+1}]}(s)ds\Bigl|^2.
\end{eqnarray}
We get from the H\"older inequality and Lemma \ref{EST_SOL}  for all $t\in [0,T]$ that
\begin{equation}
\label{A1_EST}
	\mathbb{E}|\tilde A^M_{m,1}(t)|^2 \leq TK^2\sum_{i=0}^{m-1}\int\limits_{t_i}^{t_{i+1}} (s-t_i)^2\cdot\mathbb{E}(1+|X(s)|)^2ds\leq C \max\limits_{0\leq i\leq m-1}(t_{i+1}-t_i)^2.
\end{equation}
From Lemma \ref{ITO_FORMULA} applied to $U(x)=a(t_i,x)$ and (\ref{COMP_POISS}) we have that for $s\in [t_i,t_{i+1}]$
\begin{eqnarray}
	a(t_i,X(s))-a(t_i,X(t_i))&=&\int\limits_{t_i}^s\Bigl( \frac{\partial a}{\partial y}(t_i,X(u))\cdot a(u,X(u))+\frac{1}{2}\frac{\partial^2 a}{\partial y^2}(t_i,X(u))\cdot b^2(u,X(u))\Bigr)du\notag\\
	&&+\int\limits_{t_i}^s \frac{\partial a}{\partial y}(t_i,X(u))\cdot b(u,X(u)) dW(u)\notag\\
	&&+\int\limits_{t_i}^s\Bigl(a(t_i,X(u-)+c(u,X(u-)))-a(t_i,X(u-))\Bigr)d\tilde N(u)\notag\\
	&&+\int\limits_{t_i}^s\Bigl(a(t_i,X(u-)+c(u,X(u-)))-a(t_i,X(u-))\Bigr)\lambda(u) du.\notag
\end{eqnarray}
We denote for $f\in\{a,b,c\}$ and $u \in (t_i,t_{i+1}]$
\begin{eqnarray}
	&&\alpha_i(f,u)=\frac{\partial f}{\partial y}(t_i,X(u))\cdot a(u,X(u))+\frac{1}{2}\frac{\partial^2 f}{\partial y^2}(t_i,X(u))\cdot b^2(u,X(u)), \notag\\
	&&\beta_i(f,u)=\frac{\partial f}{\partial y}(t_i,X(u))\cdot b(u,X(u)),\notag\\
	&&\gamma_i(f,u)=f(t_i,X(u-)+c(u,X(u-)))-f(t_i,X(u-)).\notag
\end{eqnarray}
We have 
\begin{equation}
\label{A2M_EST}
	\mathbb{E}|\tilde A^M_{m,2}(t)|^2\leq  4\Bigl(\mathbb{E}|\tilde M^M_{m,1}(t)|^2+\mathbb{E}|\tilde M^M_{m,2}(t)|^2+\mathbb{E}|\tilde M^M_{m,3}(t)|^2+\mathbb{E}|\tilde M^M_{m,4}(t)|^2\Bigr),
\end{equation}
where
\begin{eqnarray}
	\label{ME_1}
	&&\mathbb{E}|\tilde M^M_{m,1}(t)|^2=\mathbb{E}\Bigl|\int\limits_0^t\sum\limits_{i=0}^{m-1}\Bigl(\int\limits_{t_i}^s \alpha_i(a,u)du\Bigr)\cdot\mathbf{1}_{(t_i,t_{i+1}]}(s)ds\Bigl|^2,\\
	\label{ME_2}
	&&\mathbb{E}|\tilde M^M_{m,2}(t)|^2=\mathbb{E}\Bigl|\int\limits_0^t\sum\limits_{i=0}^{m-1}\Bigl(\int\limits_{t_i}^s \beta_i(a,u)dW(u)\Bigr)\cdot\mathbf{1}_{(t_i,t_{i+1}]}(s)ds\Bigl|^2,\\
	\label{ME_3}
	&&\mathbb{E}|\tilde M^M_{m,3}(t)|^2=\mathbb{E}\Bigl|\int\limits_0^t\sum\limits_{i=0}^{m-1}\Bigl(\int\limits_{t_i}^s\gamma_i(a,u)d\tilde N(u)\Bigr)\cdot\mathbf{1}_{(t_i,t_{i+1}]}(s)ds\Bigl|^2,\\\
	\label{ME_4}
	&&\mathbb{E}|\tilde M^M_{m,4}(t)|^2=\mathbb{E}\Bigl|\int\limits_0^t\sum\limits_{i=0}^{m-1}\Bigl(\int\limits_{t_i}^s\gamma_i(a,u)\cdot \lambda(u) du\Bigr)\cdot\mathbf{1}_{(t_i,t_{i+1}]}(s)ds\Bigl|^2,
\end{eqnarray}
for all $t\in [0,T]$. By the H\"older inequality, (\ref{LIN_GROW1}), (\ref{BOUND_AX}) and Lemma \ref{EST_SOL} we have
\begin{equation}
\label{ME_EST_1}
	\mathbb{E}|\tilde M^M_{m,1}(t)|^2\leq C_1\sum_{i=0}^{m-1}\int\limits_{t_i}^{t_{i+1}}\Bigl((s-t_i)\int\limits_{t_i}^s\mathbb{E}(1+|X(u)|)^4 du\Bigr)ds\leq C_3\max\limits_{0\leq i\leq m-1}(t_{i+1}-t_i)^2.
\end{equation}
By Theorem 6.5.8 in \cite{KUO} and Theorem 88 (iii) in \cite{situ} we obtain for all $t\in [0,T]$
\begin{eqnarray}
\label{ME_EST_2}
	\mathbb{E}|\tilde M^M_{m,2}(t)|^2
	&=&\int\limits_0^t\int\limits_0^t\sum_{i=0}^{m-1}\mathbb{E}\Bigl(\int\limits_{t_i}^{\min\{s_1,s_2\}}|\beta_i(a,u)|^2du\Bigr)\cdot\mathbf{1}_{(t_i,t_{i+1}]}(s_1)\cdot\mathbf{1}_{(t_i,t_{i+1}]}(s_2)ds_1ds_2\notag\\
&\leq& C\max\limits_{0\leq i\leq m-1}(t_{i+1}-t_i)^2,
\end{eqnarray}
\begin{eqnarray}
\label{ME_EST_23}
	\mathbb{E}|\tilde M^M_{m,3}(t)|^2
	&=&\int\limits_0^t\int\limits_0^t\sum_{i=0}^{m-1}\mathbb{E}\Bigl(\int\limits_{t_i}^{\min\{s_1,s_2\}}|\gamma_i(a,u)|^2\cdot\lambda(u)du\Bigr)\cdot\mathbf{1}_{(t_i,t_{i+1}]}(s_1)\cdot\mathbf{1}_{(t_i,t_{i+1}]}(s_2)ds_1ds_2\notag\\
 &\leq& C\max\limits_{0\leq i\leq m-1}(t_{i+1}-t_i)^2.
\end{eqnarray}
We estimate (\ref{ME_3}) analogously as (\ref{ME_1}) and we get for all $t\in [0,T]$ that
\begin{equation}
\label{ME_EST_3}
	\mathbb{E}|\tilde M^M_{m,4}(t)|^2\leq C\max\limits_{0\leq i\leq m-1}(t_{i+1}-t_i)^2.
\end{equation}
Hence, by (\ref{A2M_EST}), (\ref{ME_EST_1}), (\ref{ME_EST_2}) and (\ref{ME_EST_3}) we arrive at
\begin{equation}
	\label{A2_EST}
		\mathbb{E}|\tilde A^M_{m,2}(t)|^2\leq C\max\limits_{0\leq i\leq m-1}(t_{i+1}-t_i)^2,
\end{equation}
for all $t\in [0,T]$. For (\ref{EST_AE_3}) we have by the H\"older inequality and (A2) that
\begin{equation}
	\label{A3_EST}
	\mathbb{E}|\tilde A^M_{m,3}(t)|^2\leq \lambda^2 KT\int\limits_0^t\sum_{i=0}^{m-1}\mathbb{E}|X(t_i)-\tilde X^M_m(t_i)|^2\cdot\mathbf{1}_{(t_i,t_{i+1}]}(s)ds,
\end{equation}
for all $t\in [0,T]$. Hence, (\ref{DIFF_A_AE}), (\ref{A1_EST}), (\ref{A2_EST}) and (\ref{A3_EST}) yield for all $t\in [0,T]$ that
\begin{equation}
\label{EST_DIFF_A_AE_2}
	\mathbb{E}|A(t)-\tilde A^M_m(t)|^2\leq C_1\int\limits_0^t\sum_{i=0}^{m-1}\mathbb{E}|X(t_i)-\tilde X^M_m(t_i)|^2\cdot\mathbf{1}_{(t_i,t_{i+1}]}(s)ds+C_2\max\limits_{0\leq i\leq m-1}(t_{i+1}-t_i)^2.
\end{equation}
We have for all $t\in [0,T]$ that
\begin{equation}
	\label{DIFF_A_BM}
	\mathbb{E}|B(t)-\tilde B^M_m(t)|^2\leq 3\Bigl(\mathbb{E}|\tilde B^M_{m,1}(t)|^2+\mathbb{E}|\tilde B^M_{m,2}(t)|^2+\mathbb{E}|\tilde B^M_{m,3}(t)|^2\Bigr),
\end{equation}
where
\begin{eqnarray}
	&&\mathbb{E}|\tilde B^M_{m,1}(t)|^2=\mathbb{E}\Bigl|\int\limits_0^t\sum\limits_{i=0}^{m-1}\Bigl(b(s,X(s))-b(t_i,X(s))\Bigr)\cdot\mathbf{1}_{(t_i,t_{i+1}]}(s)dW(s)\Bigl|^2,\notag\\
	&&\mathbb{E}|\tilde B^M_{m,2}(t)|^2=\mathbb{E}\Bigl|\int\limits_0^t\sum\limits_{i=0}^{m-1}\Bigl(b(t_i,X(s))-b(t_i,X(t_i))\notag\\
	&&\quad\quad\quad\quad\quad\quad-\int\limits_{t_i}^sL_1b(U_i)dW(u)-\int\limits_{t_i}^sL_{-1}b(U_i)dN(u)\Bigr)\cdot\mathbf{1}_{(t_i,t_{i+1}]}(s)dW(s)\Bigl|^2,\notag\\
	\label{EST_BM_3}
	&&\mathbb{E}|\tilde B^M_{m,3}(t)|^2=\mathbb{E}\Bigl|\int\limits_0^t\sum\limits_{i=0}^{m-1}\Bigl(b(t_i,X(t_i))-b(U_i)\Bigr)\cdot\mathbf{1}_{(t_i,t_{i+1}]}(s)dW(s)\Bigl|^2.\notag
\end{eqnarray}
From the It\^o isometry and the H\"older inequality we obtain for $t\in [0,T]$
\begin{equation}
	\label{est_mbm_1}
\mathbb{E}|\tilde B^M_{m,1}(t)|^2\leq\sum\limits_{i=0}^{m-1}\int\limits_{t_i}^{t_{i+1}}\mathbb{E}|b(s,X(s))-b(t_i,X(s))|^2ds\leq C\max_{0\leq i\leq m-1}(t_{i+1}-t_i)^2.
\end{equation}
By the It\^o isometry together with the It\^o formula we get
\begin{equation}
\mathbb{E}|\tilde B^M_{m,2}(t)|^2\leq  4\Bigl(\mathbb{E}|\bar M^M_{m,1}(t)|^2+\mathbb{E}|\bar M^M_{m,2}(t)|^2+\mathbb{E}|\bar M^M_{m,3}(t)|^2+\mathbb{E}|\bar M^M_{m,4}(t)|^2\Bigr),
\end{equation}
where
\begin{eqnarray}
	\label{BMM_1}
	&&\mathbb{E}|\bar M^M_{m,1}(t)|^2=\mathbb{E}\int\limits_0^t\sum\limits_{i=0}^{m-1}\Bigl(\int\limits_{t_i}^s \alpha_i(b,u)du\Bigr)^2\cdot\mathbf{1}_{(t_i,t_{i+1}]}(s)ds,\\
	\label{BMM_2}
	&&\mathbb{E}|\bar M^M_{m,2}(t)|^2=\mathbb{E}\int\limits_0^t\sum\limits_{i=0}^{m-1}\Bigl(\int\limits_{t_i}^s (\beta_i(b,u)-L_1b(U_i))dW(u)\Bigr)^2\cdot\mathbf{1}_{(t_i,t_{i+1}]}(s)ds,\\
	\label{BMM_3}
	&&\mathbb{E}|\bar M^M_{m,3}(t)|^2=\mathbb{E}\int\limits_0^t\sum\limits_{i=0}^{m-1}\Bigl(\int\limits_{t_i}^s(\gamma_i(b,u)-L_{-1}b(U_i))d\tilde N(u)\Bigr)^2\cdot\mathbf{1}_{(t_i,t_{i+1}]}(s)ds,\\
		\label{BMM_4}
	&&\mathbb{E}|\bar M^M_{m,4}(t)|^2=\mathbb{E}\int\limits_0^t\sum\limits_{i=0}^{m-1}\Bigl(\int\limits_{t_i}^s(\gamma_i(b,u)-L_{-1}b(U_i))\cdot \lambda(u) du\Bigr)^2\cdot\mathbf{1}_{(t_i,t_{i+1}]}(s)ds,
\end{eqnarray}
for all $t\in [0,T]$. From the H\"older inequality we get 
\begin{equation}
	\mathbb{E}|\bar M^M_{m,1}(t)|^2\leq C(1+\sup\limits_{t\in [0,T]}\mathbb{E}|X(t)|^4)\max\limits_{0\leq i\leq m-1}(t_{i+1}-t_i)^2,
\end{equation}
for all $t\in [0,T]$. Since we have for $u\in [t_i,t_{i+1}]$ that
\begin{equation}
\mathbb{E}|\beta_i(b,u)-L_1b(U_i)|^2\leq C\Bigl((1+\sup\limits_{t\in [0,T]}\mathbb{E}|X(t)|^2)\cdot (u-t_i)^2+(u-t_i)+\mathbb{E}|X(t_i)-\tilde X^M_m(t_i)|^2\Bigr),
\end{equation}
we obtain
\begin{equation}
	\mathbb{E}|\bar M^M_{m,2}(t)|^2\leq  C_1\int\limits_0^t\sum_{i=0}^{m-1}\mathbb{E}|X(t_i)-\tilde X^M_m(t_i)|^2\cdot\mathbf{1}_{(t_i,t_{i+1}]}(s)ds+C_2\max\limits_{0\leq i\leq m-1}(t_{i+1}-t_i)^2.
\end{equation}
Moreover, for $s\in [t_i,t_{i+1}]$ we have that
\begin{equation}
	\mathbb{E}\int\limits_{t_i}^s |\gamma_i(b,u)-L_{-1}b(U_i)|^2 du \leq K_1\cdot\mathbb{E}|X(t_i)-\tilde X^M_m(t_i)|^2+K_2\cdot (t_{i+1}-t_i)^2,
\end{equation}
which implies
\begin{eqnarray}
	&&\mathbb{E}|\bar M^M_{m,3}(t)|^2= \int\limits_0^t\sum\limits_{i=0}^{m-1}\Bigl(\int\limits_{t_i}^s \mathbb{E}|\gamma_i(b,u)-L_{-1}b(U_i)|^2\cdot\lambda(u) du\Bigr)\cdot\mathbf{1}_{(t_i,t_{i+1}]}(s)ds\notag\\
&&\leq  C_1\int\limits_0^t\sum_{i=0}^{m-1}\mathbb{E}|X(t_i)-\tilde X^M_m(t_i)|^2\cdot\mathbf{1}_{(t_i,t_{i+1}]}(s)ds+C_2\max\limits_{0\leq i\leq m-1}(t_{i+1}-t_i)^2,
\end{eqnarray}
\begin{eqnarray}
	&&\mathbb{E}|\bar M^M_{m,4}(t)|^2= \int\limits_0^t\sum\limits_{i=0}^{m-1}\Bigl((s-t_i)\cdot\int\limits_{t_i}^s \mathbb{E}|\gamma_i(b,u)-L_{-1}b(U_i)|^2\cdot (\lambda(u))^2 du\Bigr)\cdot\mathbf{1}_{(t_i,t_{i+1}]}(s)ds\notag\\
&&\leq  C_1\int\limits_0^t\sum_{i=0}^{m-1}\mathbb{E}|X(t_i)-\tilde X^M_m(t_i)|^2\cdot\mathbf{1}_{(t_i,t_{i+1}]}(s)ds+C_2\max\limits_{0\leq i\leq m-1}(t_{i+1}-t_i)^3.
\end{eqnarray}
Hence, for all $t\in [0,T]$ we get
\begin{equation}
\label{est_mbm_2}
	\mathbb{E}|\tilde B^M_{m,2}(t)|^2\leq K^2 \int\limits_{0}^{t}\sum\limits_{i=0}^{m-1}\mathbb{E}|X(t_i)-\tilde X^M_m(t_i)|^2\cdot \mathbf{1}_{(t_i,t_{i+1}]}(s)ds+C_2\max\limits_{0\leq i\leq m-1}(t_{i+1}-t_i)^2,
\end{equation}
\begin{equation}
\label{est_mbm_3}
\mathbb{E}|\tilde B^M_{m,3}(t)|^2\leq K^2 \int\limits_{0}^{t}\sum\limits_{i=0}^{m-1}\mathbb{E}|X(t_i)-\tilde X^M_m(t_i)|^2\cdot \mathbf{1}_{(t_i,t_{i+1}]}(s)ds.
\end{equation}
Therefore, by (\ref{DIFF_A_BM}), (\ref{est_mbm_1}), (\ref{est_mbm_2}) and (\ref{est_mbm_3}) we obtain for all $t\in [0,T]$
\begin{equation}
\label{EST_DIFF_B_BM_2}
	\mathbb{E}|B(t)-\tilde B^M_m(t)|^2\leq C_1\int\limits_0^t\sum_{i=0}^{m-1}\mathbb{E}|X(t_i)-\tilde X^M_m(t_i)|^2\cdot\mathbf{1}_{(t_i,t_{i+1}]}(s)ds+C_2\max\limits_{0\leq i\leq m-1}(t_{i+1}-t_i)^2.
\end{equation}
Now
\begin{equation}
	\label{DIFF_A_CM}
	\mathbb{E}|C(t)-\tilde C^M_m(t)|^2\leq 3\Bigl(\mathbb{E}|\tilde C^M_{m,1}(t)|^2+\mathbb{E}|\tilde C^M_{m,2}(t)|^2+\mathbb{E}|\tilde C^M_{m,3}(t)|^2\Bigr),
\end{equation}
where
\begin{eqnarray}
	&&\mathbb{E}|\tilde C^M_{m,1}(t)|^2=\mathbb{E}\Bigl|\int\limits_0^t\sum\limits_{i=0}^{m-1}\Bigl(c(s,X(s-))-c(t_i,X(s-))\Bigr)\cdot\mathbf{1}_{(t_i,t_{i+1}]}(s)dN(s)\Bigl|^2,\notag\\
	&&\mathbb{E}|\tilde C^M_{m,2}(t)|^2=\mathbb{E}\Bigl|\int\limits_0^t\sum\limits_{i=0}^{m-1}\Bigl(c(t_i,X(s-))-c(t_i,X(t_i))\notag\\
	&&\quad\quad\quad\quad\quad\quad-\int\limits_{t_i}^s L_1c(U_i)dW(u)-\int\limits_{t_i}^{s-}L_{-1}c(U_i)dN(u)\Bigr)\cdot\mathbf{1}_{(t_i,t_{i+1}]}(s)dN(s)\Bigl|^2,\notag\\
	\label{EST_cM_3}
	&&\mathbb{E}|\tilde C^M_{m,3}(t)|^2=\mathbb{E}\Bigl|\int\limits_0^t\sum\limits_{i=0}^{m-1}\Bigl(c(t_i,X(t_i))-c(U_i)\Bigr)\cdot\mathbf{1}_{(t_i,t_{i+1}]}(s)dN(s)\Bigl|^2,\notag
\end{eqnarray}
for all $t\in [0,T]$. Next, we use the decomposition $dN(t)=d\tilde N(t)+dm(t)=d\tilde N(t)+\lambda(t)dt$ and the martingale isometry. Then the estimation  of the above terms goes in analogous way as for $\mathbb{E}|B(t)-\tilde B^M_m(t)|^2$, hence, we skip it. We get for all $t\in [0,T]$ that
\begin{equation}
\label{EST_DIFF_B_CM_2}
	\mathbb{E}|C(t)-\tilde C^M_m(t)|^2\leq C_1\int\limits_0^t\sum_{i=0}^{m-1}\mathbb{E}|X(t_i)-\tilde X^M_m(t_i)|^2\cdot\mathbf{1}_{(t_i,t_{i+1}]}(s)ds+C_2\max\limits_{0\leq i\leq m-1}(t_{i+1}-t_i)^2.
\end{equation}
Combining (\ref{EST_DIFF_A_AE_2}), (\ref{EST_DIFF_B_BM_2}) and (\ref{EST_DIFF_B_CM_2}) we get for all $t\in [0,T]$
\begin{equation}
	\sup\limits_{0\leq s\leq t}\mathbb{E}|X(s)-\tilde X^M_m(s)|^2\leq C_1\int\limits_0^t\sup\limits_{0\leq u\leq s}\mathbb{E}|X(u)-\tilde X^M_m(u)|^2 ds+C_2\max\limits_{0\leq i\leq m-1}(t_{i+1}-t_i)^2.
\end{equation}
By (\ref{EST_SOL_E}) and (\ref{finit_2mom_Mil}) the mapping $[0,T]\ni t\to \sup\limits_{0\leq s\leq t}\mathbb{E}|X(s)-\tilde X^M_m(s)|^2\in\mathbb{R}_+\cup\{0\}$ is bounded and Borel measurable. Hence, by Gronwall's lemma we get (\ref{EST_EULER_APP2}). The estimate (\ref{EST_EULER_APP1}) is a consequence of (\ref{EST_SOL_E}) and (\ref{EST_EULER_APP2}). This ends the proof. \ \ \ $\blacksquare$ 
\begin{lem} 
\label{lem_est_rm}
Let us assume that the mappings $a$, $b$, $c$ and $\lambda$ satisfy the assumptions $(A)$-$(E)$.  For all $t\in [t_i,t_{i+1}]$, $i=0,1,\ldots,m-1$
\begin{equation}
	\label{est_rmm_1}
	\mathbb{E}|\tilde R^M_m(t)|^2\leq C (t_{i+1}-t_i)^2,
\end{equation}
where $C>0$ does not depend on $m$ nor $i$.
\end{lem}
{\bf Proof.\rm}  From (\ref{ling_g_l11}) and Theorem \ref{ERR_MIL_APP} we have for $f\in\{b,c\}$ and $j\in\{-1,1\}$ that
\begin{equation}
	\label{EST_INT3}
	\mathbb{E}|L_jf(U_i)|^2\leq C,
\end{equation}
where $C>0$ does not depend on $m$ nor $i$. Moreover, for $f\in\{b,c\}$ and $j\in\{-1,1\}$ the random variable $L_jf(U_i)$ 
is $\mathcal{F}_{t_i}$-measurable. From Fact \ref{mul_int_ind} (ii) and by (\ref{COND_I_NN_START})-(\ref{COND_I_NN}) we have that $I_{t_i,t}(N,N)-\mathbb{E}(I_{t_i,t}(N,N) \ | \ \mathcal{N}_m(N))$, $I_{t_i,t}(W,W)-\mathbb{E}(I_{t_i,t}(W,W) \ | \ \mathcal{N}_m(W))$ and $I_{t_i,t}(N,W)+I_{t_i,t}(W,N)-\mathbb{E}(I_{t_i,t}(N,W)+I_{t_i,t}(W,N) \ | \ \mathcal{N}_m(N,W))$ are independent of $\mathcal{F}_{t_i}$.
Hence, by (\ref{def_RM}), Fact \ref{mul_int_ind} (i) and (\ref{EST_INT3}) we get 
\begin{eqnarray}
	\|\tilde R^M_m(t)\|_{L^2(\Omega)}&\leq & 2\cdot\|L_1b(U_i)\|_{L^2(\Omega)}\cdot\|I_{t_i,t}(W,W)\|_{L^2(\Omega)}+2\cdot\|L_{-1}c(U_i)\|_{L^2(\Omega)}\cdot\|I_{t_i,t}(N,N)\|_{L^2(\Omega)}\notag\\
	&&+2\cdot\|L_{1}c(U_i)\|_{L^2(\Omega)}\cdot\|I_{t_i,t}(N,W)+I_{t_i,t}(W,N)\|_{L^2(\Omega)}\notag\\
	&&\leq C(t-t_i),
\end{eqnarray}
for $t\in [t_i,t_{i+1}]$, which ends the proof of (\ref{est_rmm_1}). \ \ \ $\blacksquare$
\section{References}
\begingroup
\renewcommand{\section}[2]{}%

\endgroup
\end{document}